\documentclass{gtpart}
 \usepackage{amscd,amsfonts,amssymb,
amscd,amsmath,latexsym}

\title{Real secondary index theory}

\author{Ulrich Bunke}
\givenname{Ulrich}
\surname{Bunke}
\address{Universit\"at Regensburg\\ Mathematische Fakult\"at\\ 93040 Regensburg\\
Germany}
\email{ulrich.bunke@mathematik.uni-regensburg.de}
\urladdr{http://www.mathematik.uni-regensburg.de/Bunke/}

\author{Thomas Schick}
\givenname{Thomas}
\surname{Schick}
\address{Georg-August-Universit\"at G\"ottingen\\ Mathematisches Institut\\
  Bunsenstr. 3\\ 37073 G\"ottingen\\ Germany}
\email{schick@uni-math.gwdg.de}
\urladdr{http://www.uni-math.gwdg.de/schick}


\arxivreference{math.GT/0309417}
\arxivpassword{zxden}
\makeatletter

\renewcommand{\subsubsection}{\@startsection
{subsubsection}
{1}
{0mm}
{0mm}
{0mm}
{\normalfont\normalsize\itshape}}
\makeatother
\DeclareMathOperator{\im}{im}      %image

\begin{document}

% Definitionen
\newcommand{\Fred}{{\tt Fred}}
\newcommand{\Comp}{{\tt K}}
\newcommand{\reals}{\mathbb{R}}
\newcommand{\complexs}{\mathbb{C}}
\newcommand{\integers}{\mathbb{Z}}
\newcommand{\naturals}{\mathbb{N}}
\newcommand{\rationals}{\mathbb{Q}}
\newcommand{\quaternions}{\mathbb{H}}
\newcommand{\tensor}{\otimes}
\newcommand{\Proof}{{\em Proof. }}
\newcommand{\vC}{\v{C}}

\newcommand{\iso}{\cong}
\newcommand{\into}{\hookrightarrow}
\newcommand{\DD}{\mathbf{v}}
\newcommand{\orient}{{\rm or}}
\newcommand{\cSet}{\mathcal{S}et}
\newcommand{\by}{{\bf y}}
\newcommand{\bE}{{\bf E}}
\newcommand{\Face}{{\tt Face}}
\newcommand{\cDelta}{\mathbf{\Delta}}
\newcommand{\LIM}{{\tt LIM}}
\newcommand{\diag}{{\tt diag}}
 \newcommand{\dist}{{\tt dist}}
\newcommand{\kaaa}{{\frak k}}
\newcommand{\paaa}{{\frak p}}
\newcommand{\vp}{{\varphi}}
\newcommand{\taaa}{{\frak t}}
\newcommand{\haaa}{{\frak h}}
\newcommand{\Hh}{{\bf H}}
\newcommand{\Rep}{{\tt Rep}}
\newcommand{\Hb}{\mathbb{H}}
\newcommand{\str}{{\tt str}}
\newcommand{\Ind}{{\tt ind}}
\newcommand{\triv}{{\tt triv}}
\newcommand{\bD}{{\bf D}}
\newcommand{\bF}{{\bf F}}
\newcommand{\tX}{{\tt X}}
\newcommand{\Cliff}{{\tt Cliff}}
\newcommand{\tY}{{\tt Y}}
\newcommand{\tZ}{{\tt Z}}
\newcommand{\tV}{{\tt V}}
\newcommand{\tR}{{\tt R}}
\newcommand{\Fam}{{\tt Fam}}
\newcommand{\Cusp}{{\tt Cusp}}
\newcommand{\bT}{{\bf T}}
\newcommand{\bK}{{\bf K}}
\newcommand{\bo}{{\bf o}}
\newcommand{\K}{\mathbb{K}}
\newcommand{\tH}{{\tt H}}
\newcommand{\bS}{\mathbf{S}}
\newcommand{\bB}{\mathbf{B}}
\newcommand{\tW}{{\tt W}}
\newcommand{\tF}{{\tt F}}
\newcommand{\bA}{\mathbf{A}}
\newcommand{\bL}{{\bf L}}
 \newcommand{\bom}{\mathbf{\Omega}}
\newcommand{\bundle}{{\tt bundle}}
\newcommand{\ch}{\mathbf{ch}}
\newcommand{\ve}{{\varepsilon}}
\newcommand{\gen}{{\tt gen}}
\newcommand{\cTop}{\mathcal{T}op}
\newcommand{\bP}{\mathbf{P}}
\newcommand{\Naaa}{\mathbf{N}}
\newcommand{\image}{{\tt image}}
\newcommand{\gaaa}{{\frak g}}
\newcommand{\zaaa}{{\frak z}}
\newcommand{\saaa}{{\frak s}}
\newcommand{\laaa}{{\frak l}}
\newcommand{\bN}{\mathbf{N}}
\newcommand{\stimes}{{\times\hspace{-1mm}\bf |}}
\newcommand{\ausg}{{\rm end}}
\newcommand{\bff}{{\bf f}}
\newcommand{\maaa}{{\frak m}}
\newcommand{\aaaa}{{\frak a}}
\newcommand{\naaa}{{\frak n}}
\newcommand{\brr}{{\bf r}}
\newcommand{\res}{{\tt res}}
\newcommand{\Aut}{{\tt Aut}}
\newcommand{\Pol}{{\tt Pol}}
\newcommand{\Tr}{{\tt Tr}}
\newcommand{\cT}{\mathcal{T}}
\newcommand{\dom}{{\tt dom}}
\newcommand{\Line}{{\tt Line}}
\newcommand{\db}{{\bar{\partial}}}
\newcommand{\Sf}{{\tt  Sf}}
\newcommand{\g}{{\gaaa}}
\newcommand{\cZ}{\mathcal{Z}}
\newcommand{\cH}{\mathcal{H}}
\newcommand{\cM}{\mathcal{M}}
\newcommand{\interi}{{\ttint}}
\newcommand{\singsupp}{{\tt singsupp}}
\newcommand{\cE}{\mathcal{E}}
\newcommand{\ccR}{\mathcal{R}}
\newcommand{\hol}{{\tt hol}}
\newcommand{\cV}{\mathcal{V}}
\newcommand{\cY}{\mathcal{Y}}
\newcommand{\cW}{\mathcal{W}}
\newcommand{\dR}{{\tt dR}}
\newcommand{\del}{{\tt del}}
\newcommand{\bdel}{\mathbf{del}}
\newcommand{\cI}{\mathcal{I}}
\newcommand{\cC}{\mathcal{C}}
\newcommand{\cK}{\mathcal{K}}
\newcommand{\cA}{\mathcal{A}}
\newcommand{\cU}{\mathcal{U}}
\newcommand{\Hom}{{\tt Hom}}
\newcommand{\vol}{{\tt vol}}
\newcommand{\cO}{\mathcal{O}}
\newcommand{\End}{{\tt End}}
\newcommand{\Ext}{{\tt Ext}}
\newcommand{\sign}{{\tt sign}}
\newcommand{\spann}{{\tt span}}
\newcommand{\symm}{{\tt symm}}
\newcommand{\cF}{\mathcal{F}}
\newcommand{\cD}{\mathcal{D}}
\newcommand{\bC}{\mathbf{C}}
\newcommand{\bbeta}{\mathbf{\eta}}
\newcommand{\bOmega}{\mathbf{\Omega}}
\newcommand{\bbz}{{\bf z}}
\newcommand{\bc}{\mathbf{c}}
\newcommand{\bb}{\mathbf{b}}
\newcommand{\bd}{\mathbf{d}}
\newcommand{\Ree}{{\tt Re }}
\newcommand{\Res}{{\tt Res}}
\newcommand{\Imm}{{\tt Im}}
\newcommand{\inter}{{\tt int}}
\newcommand{\clo}{{\tt clo}}
\newcommand{\tg}{{\tt tg}}
\newcommand{\ee}{{\tt e}}
\newcommand{\Li}{{\tt Li}}
\newcommand{\cN}{\mathcal{N}}
 \newcommand{\conv}{{\tt conv}}
\newcommand{\op}{{\tt Op}}
\newcommand{\tr}{{\tt tr}}
\newcommand{\ctg}{{\tt ctg}}
\newcommand{\degg}{{\tt deg}}
\newcommand{\Ad}{{\tt  Ad}}
\newcommand{\ad}{{\tt ad}}
\newcommand{\codim}{{\tt codim}}
\newcommand{\Gr}{{\tt Gr}}
\newcommand{\coker}{{\tt coker}}
\newcommand{\id}{{\tt id}}
\newcommand{\ord}{{\tt ord}}
\newcommand{\nat}{\mathbb{N}}
\newcommand{\supp}{{\tt supp}}
\newcommand{\sing}{{\tt sing}}
\newcommand{\spec}{{\tt spec}}
\newcommand{\Ann}{{\tt Ann}}
 \newcommand{\Or}{{\tt Or }}
\newcommand{\Diff}{\mathcal{D}iff}
\newcommand{\cB}{\mathcal{B}}
\def\imath{{i}}
\newcommand{\cR}{\mathcal{R}}
\def\hB{\hspace*{\fill}$\Box$ \newline\noindent}
\newcommand{\varho}{\varrho}
\newcommand{\ind}{{\tt index}}
\newcommand{\Indu}{{\tt Ind}}
\newcommand{\Fin}{{\tt Fin}}
\newcommand{\cS}{\mathcal{S}}
\newcommand{\orig}{\mathcal{O}}
\def\hB{\hspace*{\fill}$\Box$ \\[0.5cm]\noindent}
\newcommand{\cL}{\mathcal{L}}
 \newcommand{\cG}{\mathcal{G}}
\newcommand{\Mat}{{\TT Mat}}
%\newcommand{\npi}{{\naaa_P\hspace{-1.5mm}-\hspace{-2mm}\mbox{\rm inv}}}
%\newcommand{\ngp}{{N_\Gamma(\pi)}}
%\newcommand{\gbg}{{\Gamma\backslash G}}
%\newcommand{\gkm}{{ Mod(\gaaa,K) }}
%\newcommand{\ggkm}{{  (\gaaa,K) }}
%\newcommand{\pkm}{{ Mod(\paaa,K_P)}}
%newcommand{\ppkm}{{  (\paaa,K_P)}}
%\newcommand{\makm}{{Mod(\maaa_P\oplus\aaaa_P,K_P)}}
%\newcommand{\mmakm}{{ (\maaa_P\oplus\aaaa_P,K_P)}}
\newcommand{\cP}{\mathcal{P}}
\newcommand{\bv}{\mathbf{v}}
\newcommand{\cQ}{\mathcal{Q}}
 \newcommand{\cX}{\mathcal{X}}
\newcommand{\bH}{\mathbf{H}}
\newcommand{\bW}{\mathbf{W}}
\newcommand{\pr}{{\tt pr}}
\newcommand{\bX}{\mathbf{X}}
\newcommand{\bY}{\mathbf{Y}}
\newcommand{\bZ}{\mathbf{Z}}
\newcommand{\bz}{\mathbf{z}}
\newcommand{\bkappa}{\mathbf{\kappa}}
\newcommand{\ev}{{\tt ev}}
\newcommand{\bV}{\mathbf{V}}
\newcommand{\Gerbe}{{\tt Gerbe}}
\newcommand{\gerbe}{{\tt gerbe}}
\newcommand{\hA}{\mathbf{\hat A}}
\newcommand{\Sets}{\mathcal{S}ets}

\newtheorem{prop}{Proposition}[section]
\newtheorem{proposition}[prop]{Proposition}
\newtheorem{lem}[prop]{Lemma}
\newtheorem{lemma}[prop]{Lemma}
\newtheorem{theorem}[prop]{Theorem}
\newtheorem{fact}[prop]{Fact}

\newtheorem{notation}[prop]{Notation}
\newtheorem{kor}[prop]{Corollary}
\newtheorem{corollary}[prop]{Corollary}
\newtheorem{ddd}[prop]{Definition}
\newtheorem{definition}[prop]{Definition}
\newtheorem{ass}[prop]{Assumption}
\newtheorem{con}[prop]{Conjecture}
\newtheorem{prob}[prop]{Problem}
\newtheorem{remark}[prop]{Remark}
\newtheorem{example}[prop]{Example}

{\catcode`@=11\global\let\c@equation=\c@prop}
\renewcommand{\theequation}{\theprop}
% Hier werden Gleichungen und Theoreme zusammen gezaehlt. Soll ein anderer Zaehler statt theorem verwendet werden (entspr. dem \newtheorem-Befehl), muss 2-mal theorem durch diesen Zaehler ersetzt werden. (Die Zeilen entsprechen der Zaehlung von \newtheorem{equation}[theorem]).

\begin{abstract}
  In this paper, we study the family index of a family of spin
  manifolds. In particular, we discuss to what extend the real index
  (of the Dirac operator of the real spinor bundle if the fiber
  dimension is divisible by $8$)
  which can be defined in this case contains extra information over
  the complex index (the index of its complexification). We study this
  question under the additional assumption that the complex index
  vanishes on the $k$-skeleton of $B$. In this case, we define new
  analytical invariants $\hat c_k\in H^{k-1}(B;\reals/\integers)$,
  i.e.~a certain secondary invariant.

  We give interesting non-trivial examples.
 We then describe this invariant in terms of known topological
  characteristic classes.
\end{abstract}

\maketitle

\tableofcontents
\parskip3ex

\section{Introduction}
\subsubsection{}
The index of a family of Fredholm operators parametrized by a space
$B$ is an element in the K-theory $K^*(B)$ of this parameter space. If the
base is in fact a smooth compact manifold without boundary, and
this family is a 
family of fiberwise generalized Dirac operators on a smooth fiber
bundle over $B$, then after adding some further geometric structures in order
to define the Bismut super connection 
we can do local index theory in the sense of
\cite{berlinegetzlervergne92}. Let us denote by $\cE$ the family with this collection
of geometric structures, by $D(\cE)$ the family of Dirac operators, and
by $\ind(\cE)$ the index of this family.
Local index theory provides a closed
form $\Omega(\cE)$ on $B$ (see Definition \ref{lofo}) which represents a cohomology class
$[\Omega(\cE)]\in H^*(B,\R)$. The local index theorem states that
$$\ch^\R(\ind(\cE))=[\Omega(\cE)]\ .$$
\subsubsection{}
The focus of the present paper is not a generalization of this type of
result. Let us illustrate the philosophy of the present paper in
the case above. We start with local index theory and produce the
even form $\Omega(\cE)$. We then observe that this form is closed and
therefore represents a cohomology class $[\Omega(\cE)]\in H^{even}(B,\R)$.
We observe that this class in fact only depends on
$\ind(\cE)\in K^0(B)$. The classifying space of  the functor
$K^0$ is $BU\times \Z$.
By naturality we conclude that there must
be a universal class $\ch^\R_{univ}\in H^{ev}(BU\times \Z,\R)$
such that $[\Omega(\cE)]=f^*\ch^\R_{univ}$, if $f:B\rightarrow BU\times
\Z$ classifies $\ind(\cE)$.
We know that $H^{ev}(BU,\R)$
is a polynomial ring in generators $c_2^\R,c_4^\R,\dots$.
Then we finally look for a formula which expresses 
$\ch^\R_{univ}$ in terms of these generators, this way obtaining the
Chern character. Of course this is a well known possible way toward
the family index formula.

\subsubsection{}
In fact the real subject of the present paper is a secondary version
of this approach. Let $K_k^*(B)$ denote the $k$'th step of the
Atiyah-Hirzebruch filtration of $K$-theory (see \ref{ahfi}), i.e.~the
subgroup of classes which vanish when restricted to the $(k-1)$-skeleton of
$B$ (so that $K_0^*(B)=K(B)$).
 Under the assumption that $\ind(\cE)\in K_{k}^*(B)$, the Chern
 class $c_k(\ind(\cE))\in H^k(B,\Z)$ 
(note that we use a non-standard notation where the subscript is equal to the degree)
admits a natural lift to smooth Deligne cohomology
$\hat c_k(\cE)\in H^k_{Del}(B)$ (see \ref{deefg1} and \ref{co2}).
This lift is a differential-geometric (or even global-analytic)
invariant which varies continuously with the geometry.
In particular it has a curvature
$\omega(\hat c_k(\cE))\in\cA^k(B)$, which can be expressed through $\Omega(\cE)$.

\subsubsection{}
We rigidify the situation be imposing additional geometric
constraints. We in fact assume that the family of Dirac operators $D(\cE)$
is a family of twisted Dirac operators on a family of
$Spin$-manifolds, and that the twisting bundle is a real bundle.
Let $n$ be the fiber-dimension of this family. It follows from the
presence of the real structure that, if $k+n\equiv 2(4)$ then the class $\hat
c_k(\cE)$ is flat (see \ref{dhg1}). This means that $\omega(\hat c_k(\cE))=0$.
Since any two geometric structures can be connected by a path we can
now conclude that under this assumptions $\hat c_k(\cE)$ is a
differential-topological invariant. In Subsection \ref{trio} we give
some non-trivial examples.

\subsubsection{}
Note that the flat part of $H^k_{Del}(B)$ can be identified with
$H^{k-1}(B,\R/\Z)$. Thus, given a family of $n$-dimensional spin
manifolds and a real twisting bundle such that $\ind(\cE)\in
K_{k}^*(B)$ and $k+n\equiv 2 \pmod 4$ we have defined a class $\hat c_k(\cE)\in
H^{k-1}(B,\R/\Z)$. This class is natural under pull-back.

The index of the family $D(\cE)$ has a real refinement
$\ind_\R(\cE)\in KO^{-n}(B)$ (see Subsection \ref{trara}). We prove that
in fact $\hat c_k(\cE)$ only depends on $\ind_\R(\cE)$
and that $2\hat c_k(\cE)=0$ (see Propositions \ref{plam} and \ref{zmal0}).

Let $U^{-n}_{l}(B)\subset KO^{-n}(B)$
be the subset of classes which after complexification belong to
$K_{l}^{-n}(B)$. What we have constructed so far is a natural
transformation 
\begin{equation*}
d^{-n}_{B,4k+1-n}:U^{-n}_{4k+2-n}(B)\rightarrow
H^{4k+1-n}(B,\R/\Z)
\end{equation*}
such that $d^{-n}_{B,4k+1-n}(\ind_\R(\cE))=\hat c_{4k+2}(\cE)$.

\subsubsection{}
The universal situation is given by the fiber sequence
$$\Omega^nU/O\stackrel{\Omega^n i}{\rightarrow} \Omega^n(BO\times
\Z)\rightarrow \Omega^n(BU\times \Z)\ .$$
It is obtained by application of the functor $\Omega^n$ to the
fibration $U/O\xrightarrow{i} BO\to BU$, where we construct $BO=EU/O$
and thus obtain the inclusion $i\colon U/O\to BO$.
A class $x\in KO^{-n}(B)$ is represented by a map
$f:B\rightarrow \Omega^n(BO\times \Z)$. If $x$ belongs to
$U^{-n}_\infty(B)$, then this map factors up to homotopy through a map
$g:B\rightarrow \Omega^n
U/O $.
Thus there must be a universal class $$\bar d_{4k+1-n}\in H^{4k+1-n}(\Omega^n
U/O,\R/\Z)$$ such that $g^*\bar d_{4k+1-n}=d^{-n}_{B,4k+1-n}(x)$.
Note that this universal class has the special property that
$g^*\bar d_{4k+1-n}$ only depends on the homotopy class of the composition
$\Omega^ni\circ g$.

\subsubsection{}
The main purpose of the present paper is 
the calculation of the universal class $\bar d_{4k+1-n}$
in terms of the classically known generators of the cohomology of
$\Omega^n U/O$.
The result is presented in Theorem \ref{theo:final answaer}.

\subsubsection{}
If $n=2,3,4,5$ or $n=1,k>0$, then $\bar d_{4k+1-n}=0$.
If $n=0$ or $n=1,k=0$ then the  class 
$d^{-n}_{B,4k+1-n}(x)$
is a classical characteristic class of $x$
(i.e.~it can be expressed in terms of the dimension and Stiefel-Whitney
classes). If $n=7$, in principle it is also a well known characteristic class
(albeit for $KO^7$ which is not considered much), i.e.~pulls back from the
classifying space $U/O$.

The dimension $n=6$ is interesting since in this case 
the classes are definitely not just classical characteristic classes of
$x$.

See Subsection \ref{cl7} for all that.

\subsubsection{}
For the convenience of the reader we have added two appendixes.
In the first we recall (with proof) some material about
transgression.

In the second we recall the results of Cartan
\cite{cartan} about the cohomology of the
spaces $\Omega^n U/O$ and about the action of various maps and 
transgressions on this cohomology.

\section{Real local index theory}

\subsection{Chern classes of geometric families in Deligne cohomology}

\subsubsection{}
We consider a fiber bundle $\pi:E\rightarrow B$ with closed $n$-dimensional fibers.
We assume that the vertical $T^v\pi:=\ker(d\pi)$ bundle is oriented
and equipped with a spin structure. We choose a vertical  Riemannian
metric $g^{T^v\pi}$ and a horizontal distribution $T^h\pi$, i.e. a
complement of $T^v\pi$ in $TE$.
Finally, we  let $\bW:=(W,\nabla^W,h^W)$ be an auxiliary complex
vector bundle
with hermitian metric and metric connection.
The data described so far make up a geometric family $\cE$ over $B$.

\subsubsection{}

We assume that $\bW$ admits a real structure $Q\in \End(W_{|\R})$
which is compatible with the connection and the metric. Then $\bW$ is the complexification of a real bundle $\bW_\R=(W_\R,\nabla^{W_\R},h^{W_\R})$. The latter can be identified with the $+1$-eigenbundle of $Q$.

\subsubsection{}

The data which we compressed in the notion of a geometric family
induces a family of elliptic operators $D(\cE)$ over $B$. Indeed, for $b\in B$
the operator $D(\cE)(b)$ is the spin Dirac operator of the Riemannian
spin manifold $E_b:=\pi^{-1}(b)$ twisted by the bundle $\bW_{|E_b}$.
The family index of $D(\cE)$ is the element
$\ind(\cE)\in K^{- n}(B)$.

\subsubsection{}

For $k\in \nat_0$ we introduce a natural transformation
of $c_k:K^{ n}(B)\rightarrow H^k(B,\Z)$ given by the Chern class.
In order to have a uniform notation in the even and odd dimensional
case we use a notation which differs from the conventional one.
So if $n$ is even, then we set
$c_{2k}:=\bc_k$ and $c_{2k+1}:=0$, where $\bc_k:K^0(B)\rightarrow
H^{2k}(B,\Z)$ is the Chern class in the usual notation.
If $n$ is odd, then we set $c_{2k}=0$ and define
$c_{2k+1}$ such that the following diagram is commutative
$$\begin{array}{ccc}
K^{ n}(B)&\stackrel{c_{2k+1}}{\rightarrow} &H^{2k+1}(B,\Z)\\
\|&& \|\\
\tilde K^{{ n+1}}(\Sigma B)&\stackrel{\bc_{2k+2}}{\rightarrow}
&H^{2k+2}(\Sigma B,\Z)
\end{array}\ ,$$
where $\tilde K$ denotes the reduced $K$-theory and the vertical
isomorphisms are the natural suspension isomorphisms.

\subsubsection{}
The Chern character is a natural transformation
$$\ch\colon K^{ n}(B)\rightarrow \bigoplus_{k\equiv n (2)}H^{ k}(B,\Q)\ .$$
Here again for the even part
$\ch:K^{2n}(B)\rightarrow  \bigoplus_{k\equiv 0({2})}H^{ k}(B,\Q)$
we use the usual convention, while the odd part is defined such that
the following diagram is commutative
$$\begin{array}{ccc}
K^{ 2n-1}(B)&\stackrel{\ch}{\rightarrow} & \bigoplus_{k\equiv 1({2})}H^{k}(B,\Q)\\
\|&& \|\\
\tilde K^{ 2n}(\Sigma B)&\stackrel{\ch}{\rightarrow}
& \bigoplus_{k\equiv 0\pmod{2}}\tilde H^{k}(\Sigma B,\Q)
\end{array}\ .$$
For $k\in\nat_0$ let
$\ch_k:K^{*}(B)\rightarrow H^k(B,\Q)$ denote the corresponding component.
%Furthermore, we extend $\ch_k$ by zero to all of $K^*(B)$.

\subsubsection{}\label{ahfi}
The ring $K^*(B)$ has a natural decreasing  filtration, the Atiyah-Hirzebruch
filtration \cite{atiyahhirzebruch},
$$\dots \subset K^*_{k+1}(B)\subset K^*_k(B)\subset\dots \subset K^*_0(B)=K^*(B)\ .$$
Recall that
$x\in K_k^*(B)$ iff $f^*x=0$ for any $k-1$-dimensional $CW$-complex
$X$ and continuous map $f:X\rightarrow B$.

\subsubsection{}

 %Let $k\in\nat_0$.
%Let $\Omega^k(\cE)$ denotes the degree $k$-component of the local index form
%of the family. If $n$ is even, then $\ind(\cE)\in K^0(B)$ and
%$c_k(\ind(\cE))$ is the Chern class of the index of the family of
%Dirac operators associated to the family $\cE$.
%This class is zero by definition for odd $k$.
%If $n$ is odd, then $\ind(\cE)\in K^{-1}(B)$.
%By definition $K^{-1}(B)=\tilde K^0(\Sigma B)$, where
%$\Sigma$ denotes the suspension of $B$ and $\tilde K^0$ is the reduced $K$-theory.
%Then $c_{k}(\ind(\cE))\in H^k(B,\Z)$ is the class which corresponds
%to $c_{k+1}(\ind(\cE))$ under the suspension isomorphism
%$H^{k+1}(\Sigma B,\Z)\cong H^k(\Sigma,\Z)$.
%These classes will be called the odd Chern classes.
%For even $k$ the odd Chern classes are zero by definition.

Fix now $k\in\nat_0$ and define $m\in\nat$ such that $k=m$ or
$k=2m-1$. If $x\in K_k^*(B)$, then we have
\begin{equation}\label{revi1}c_k(x)_\Q=(-1)^{m-1}(m-1)! \ch_k(x)\end{equation}
where $c_k(x)_\Q\in H^k(B,\Q)$ is the natural image of $c_k(x)$
in cohomology with rational coefficients.

\subsubsection{}\label{deefg1}

Let $H^*_{Del}(B)$ denote the smooth Deligne cohomology of $B$.
In the present paper we use its description in terms of differential
characters given by Cheeger-Simons \cite{cheegersimons83}. 
Let $Z^{k-1}$ be the group of smooth singular chains on $B$.
A class $\hat x\in H_{Del}^k(B)$ is a homomorphism
$\hat x:Z^{k-1}\rightarrow \R/\Z$ such that there exists a
differential form
$\omega(\hat x)\in \cA^k(B)$ with the property that for any smooth
singular $k$-chain $C$ we have $\hat x(\partial C)=[\int_C\omega(\hat
x)]$,
where $[r]\in \R/\Z$ denotes the class of $r\in \R$.
Note that $\omega(\hat x)$ is uniquely determined by $\hat x$. It is
called the curvature of $\hat x$. It is necessarily closed and has
integral periods.

The association $B\mapsto H^*_{Del}(B)$ is a contravariant functor
from smooth manifolds and smooth maps to graded abelian groups.
There is a natural exact sequence
$$H^{k-1}(B,\Z)\rightarrow \cA^{k-1}(B)/\im(d)\stackrel{\hat a}{\rightarrow}
H^k_{Del}(B)\stackrel{v}{\rightarrow} H^k(B,\Z)\rightarrow 0\ ,$$
where $a$ is given by
$$\hat a(\beta)(Z)=[\int_{Z}\alpha],\quad \beta\in \cA^{k-1}(B)\ .$$
Note that $\omega(\hat a(\beta))=d\beta$.
The map $v$ has the following description. Let $\hat x\in H^{k}_{Del}(B)$.
We choose a smooth $\R$-valued $(k-1)$-cochain $T$ such that
$T_{|Z^{k-1}}=\hat x$. This is possible since $\R$ is divisible.
Then we have $dT=\omega-c$ for some $\Z$-valued $k$-cochain $c$.
It follows that $c$ is closed, and we set $v(\hat x):=[c]$.
For details we refer to Cheeger-Simons  \cite{cheegersimons83}.

\subsubsection{}

A complex vector bundle $W\rightarrow B$ represents an element $[W]\in
K^0(B)$. Assume that $W$ comes with a hermitian metric $h^W$ and
metric connection $\nabla^W$. We set $\bW:=(W,h^W,\nabla^W)$. For
$k\in\nat_0$, Cheeger-Simons \cite{cheegersimons83} constructed a natural lift $\hat
c_{2k}(\bW)\in H^{2k}_{Del}(B)$ of
$c_{2k}([W])$ such that $v(\hat
c_{2k}(\bW))=c_{2k}([W])$ and $\omega(\hat c_{2k}(\bW))\in\cA^{2k}(B)$
is the Chern-Weyl representative of $c_{2k}([W])_\R$ associated to the
connection $\nabla^W$.

\subsubsection{}
The bundle $\bW$ can be considered as a geometric family $\cW$ over $B$ with
zero-dimensional fiber in a natural way. In this case we have $\ind(\cW)=[W]$.

Therefore we can consider the geometric family $\cE$ over $B$ as a
generalization of a hermitian vector bundle with connection over $B$.
It is now an obvious  question whether one can define a natural lift
$\hat c_{k}(\cE)\in H^k_{Del}(B)$ of $c_{k}(\ind(\cE))$.

\subsubsection{}
The geometric data associated with the geometric family $\cE$ induce a
connection $\nabla^{T^v\pi}$ on the vertical bundle in a natural
way. In fact, if we choose for a moment a Riemannian metric $g^{TB}$ on the
base, then we can define a Riemannian metric $g^{TE}$ on the total space $E$ as the
orthogonal sum of the vertical metric $g^{T^v\pi}$ and the metric
$g^{T^h\pi}$ on the horizontal
bundle obtained by lifting $g^{TB}$. Then $\nabla^{T^v\pi}$ is the
projection of the Levi-Civita connection of $g^{TE}$ to the vertical
bundle. This connection does not depend on the choice of $g^{TB}$.
We refer to \cite{berlinegetzlervergne92} for details.
By $\hA(\nabla^{T^v\pi})\in\cA^*(E)$ we denote the Chern-Weil
representative of the $\hA$-class of $T^v\pi$.
Furthermore, let $\ch(\nabla^{W})\in\cA(E)$ be the Chern-Weil
representative of the Chern character of $W$.
\begin{ddd}\label{lofo}
The local index form $\Omega(\cE)\in \cA(B)$ of the geometric family $\cE$ is defined by
$$\Omega(\cE):=\int_{E/B}\hA(\nabla^{T^v\pi})\ch(\nabla^W)\ .$$
\end{ddd}

\subsubsection{}
The Atiyah-Singer index theorem for families states that
$$\ch(\ind(\cE))_\R=[\Omega(\cE)]\ ,$$
where $[\omega]\in H^*(B,\R)$ denotes the class represented by
the closed form  $\omega$.
Here we once and for all fix the isomorphism between de Rham
cohomology $H^*_{dR}(M)$ and singular cohomology $H^*(B,\R)$
which is induced by the integration map. This means that the value of
the class $[\omega]$ on the cycle $Z$ is given by $\int_Z\omega$.

\subsubsection{}
The form $\Omega(\cE)$ plays the role of the Chern-Weil representative
of the Chern character of an index bundle with connection $\nabla^{\ind(\cE)}$ for $D(\cE)$ though we are not able
to define the latter object. In particular, the local index form also
determines candidates for the Chern-Weil representatives
$c_{k}(\nabla^{\ind(\cE)})$ of the Chern classes $c_k(\ind(\cE))$.
Unfortunately we are not able to define natural lifts
$\hat c_{k}(\cE)\in H^k_{Del}(B)$  of $c_{k}(\ind(\cE))$ with curvature
$\omega(\hat c_{k}(\cE))=c_{k}(\nabla^{\ind(\cE)})$.

\subsubsection{}
Assume that $\ind(\cE)\in K_k^*(B)$ and $k=2m$ or $k=2m-1$. 
By equation (\ref{revi1}) we have
$$c_{k}(\ind(\cE))_\R=[(-1)^{m-1}(m-1)!\Omega^{k}(\cE)]\ .$$
In \cite[Definition 8.19]{bunke020} we have constructed a natural class
$$
\hat
c_{k}(\cE)\in H_{Del}^k(B)
$$
with curvature $\omega(\hat c_k(\cE))=\Omega^k(\cE)$ and $v(\hat
c_{k}(\cE))=
c_{k}(\ind(\cE))$.
Instead of repeating the rather indirect construction \cite{bunke020}
we give here a direct description which could be taken as definition
of $\hat
c_{k}(\cE)$ as well.
Note that $\hat c_k(\cE)=0$ by definition if $n\not\equiv k(2)$
(recall that $n$ is the dimension of the fiber of $\cE$).
\subsubsection{}\label{thik}
Let $Z\in Z^{k-1}$ be a smooth cycle.
We must prescribe
$\hat c_k(\cE)(Z)$.
We can find a smooth manifold $X$ (not necessarily closed) of the homotopy type
of a $k-1$-dimensional $CW$-complex, a map $f:X\rightarrow B$, and  a
smooth $k-1$-cycle $Z^\prime$ in $X$, such that $f_*Z^\prime=Z$.
We could e.g.~take for $f:X\rightarrow B$ the inclusion of a  thickening of the
trace $|Z|\subset B$ of $Z$ and $Z=Z^\prime$.
Note that $0=f^*\ind(\cE)=\ind(f^*\cE)$. Therefore we can find a perturbation of
the family of Dirac operators $D(f^*\cE)$ by a family of selfadjoint smoothing
operators $Q$ (which are odd in the even-dimensional case) such that
the family $D(f^*\cE)+Q$ is invertible. In \cite{bunke020} the pair
$(f^*\cE,Q)$ was called a tamed geometric family and denoted by $f^*\cE_t$.

\subsubsection{}
If $\cF_t$ is a tamed geometric family over some base $B$, then the super connection formalism provides
an eta-form $\eta(\cF_t)\in\cA(B)$ such that $d\eta(\cF_t)=\Omega(\cE)$.
We refer to \cite{bunke020} and \cite{berlinegetzlervergne92} for details.
The form $\eta(\cF_t)$ depends on the taming. Assume that
$\cF_t^\prime$ is a second  taming of the same underlying geometric
family. Then the difference $\eta(\cF_t)-\eta(\cF_t^\prime)$ 
is a closed form. As a consequence of the index theorem for boundary
tamed families \cite{bunke020} we know that
$$[\eta(\cF_t)-\eta(\cF_t^\prime)]=\ch(x)_\R$$
for some $x\in K^*(B)$. In fact, we can take $x=\ind((\cF\times I)_{bt})$,
where  the boundary  taming is induced by $\cF_t$ and $\cF_t^\prime$.

\subsubsection{}\label{co2}
We can now prescribe $\hat c_k(\cE)(Z)$ as follows.
\begin{ddd}
\begin{equation}\label{eq:compute_char_class}
\hat
c_k(\cE)(Z):=[(-1)^{m-1}(m-1)!\int_{Z^\prime}\eta^{k-1}(f^*\cE_t)] \in
\R/\Z\ .\end{equation}
\end{ddd}
In order to see that  $\hat c_k(\cE)$ is well-defined note that
$(m-1)!\ch_k$ is an integral cohomology class. Therefore the right-hand side
does not depend on the choice of the taming.
One also checks independence of $f,X,Z^\prime$.
The relation
$\omega(\hat c_k(\cE))=\Omega^k(\cE)$
follows from $d\eta(f^*\cE_t)=f^*\Omega(\cE)$.

\subsubsection{}\label{req1}
Up to this point we have not employed the fact that the geometric bundle
$\bW=(W,\nabla^W,h^W)$ comes with a real structure $Q$.
Because of the existence of $Q$ the geometric bundle $\bW$ is isomorphic
to its hermitian conjugate $\bar \bW$.
We conclude from the general equality
$$\ch_{2k}(\nabla^{\bar V})=(-1)^k\ch_{2k}(\nabla^{V})$$
that  $\ch_{l}(\nabla^W)=0$ if $l\not\equiv 0(4)$.
 
\subsubsection{}\label{dhg1}
Recall that $n=\dim(E)-\dim(B)$.
\begin{lem}\label{vancur}
If $k+n\not\equiv 0(4)$, then
$\Omega^k(\cE)=0$.
\end{lem}
\begin{proof}
We have
$$\Omega(\cE)=\int_{E/B}\hA(\nabla^{T^v\pi})\ch(\nabla^W)\ .$$
Since the form $\hA(\nabla^{T^v\pi})$ is non-trivial only in degrees
$4l$, $l\ge
0$, we immediately see that
$\Omega^k(\cE)=0$ if $k+n\not\equiv 0(4)$.
\end{proof}

\subsubsection{}\label{firstbock}
We call a class $\hat x\in H^k_{Del}(B)$ flat if $\omega(\hat x)=0$.
By Lemma \ref{vancur} the class
$\hat c_{k}(\cE)$ is flat if $k+n\equiv 2(4)$.
The Deligne cohomology of $B$ fits into the exact sequence (see \cite{cheegersimons83})
\begin{equation}
0\rightarrow H^{k-1}(B,\R/\Z)\stackrel{\hat b}{\rightarrow}
H^k_{Del}(B)\stackrel{\omega}{\rightarrow}
\cA^k(B)
\label{eq:flat_sequence}
\end{equation}
such that $v\circ \hat b:H^{k-1}(B,\R/\Z)\rightarrow H^k(B,\Z)$
is the Bockstein operator associated to the exact sequence of
coefficients
$$0\rightarrow \Z\rightarrow \R\rightarrow \R/\Z\rightarrow 0\ .$$
By \eqref{eq:flat_sequence}, a flat class in $H^k_{Del}(B)$ can be considered as a class in
$H^{k-1}(B,\R/\Z)$.
From now on we consider $H^{k-1}(B,\R/\Z)$ as a subset of $H^k_{Del}(B)$
and do not write $\hat b$ anymore.

\subsubsection{}
The first assertion of the following proposition is just the
conclusion of the preceeding discussion.
\begin{prop}
\begin{enumerate}
\item
Let $\cE$ be a geometric family over $B$ such that the geometric twisting bundle
$\bW$ admits a real structure. Let $k\ge 0$ and assume that
$\ind(\cE)\in K_k^*(B)$. If furthermore $k+n\equiv2(4)$ (where
$n=\dim(E)-\dim(B)$ is the fiber dimension of $\cE$)
then $\hat c_k(\cE)\in H^k_{Del}$ is flat and therefore gives
rise to a class in
$H^{k-1}(B,\R/\Z)$. 
\item
The class $\hat c_k(\cE)\in H^{k-1}(B,\R/\Z)$
is independent of the geometric structures, i.e. it only depends on
the differentiable fiber bundle $E\rightarrow B$, the choice of spin structure
and orientation of the vertical bundle $T^v\pi$, and on the real vector bundle
$W_R$.
\end{enumerate}
\end{prop}
\begin{proof}
In order to show the independence of the geometric structures we argue using the connectedness of the space $P$ of these structures.
We can set up a universal family $\cE_{univ}$ over $P\times B$ and define
$\hat c_k(\cE_{univ})\in H^{k-1}(P\times B,\R/\Z)$. It follows from the homotopy invariance of the cohomology functor and the naturality of the construction of these classes with respect to pull-back, that the specializations of $\hat c_k(\cE_{univ})$ at different parameter
points $p,q$ are cohomologous on the one hand, and give the
classes $\hat c_k(\cE_p)$ and $\hat c_k(\cE_q)$ associated to the families
$\cE_p$ and $\cE_q$ equipped with geometric structures given by $p$
and $q$, respectively, on the other hand. \end{proof}

\subsubsection{}
The main goal of the present paper is to understand the nature
of the class $\hat c_{k}(\cE)\in H^{k-1}(B,\R/\Z)$ in terms of the
topology of the geometric family.

\subsection{The real index}\label{trara}

\subsubsection{}\label{cbintro}
The group $KO^0(B)$ is defined as the group completion of the
semigroup of isomorphism classes of real vector bundles over $B$.
The functor $B\mapsto KO^0(B)$ extends to an $8$-periodic multiplicative cohomology theory $KO^*$.
Complexification of real vector bundles induces a natural
transformation $c_B:KO^0(B)\rightarrow K^0(B)$ which extends to a
natural transformation $c_B:KO^*(B)\rightarrow KO^*(B)$
of multiplicative cohomology theories.

If $k+n\not\equiv 0(4)$, then the composition
$$KO^{-n}(B)\stackrel{c_B}{\rightarrow}
K^{-n}(B)\stackrel{\ch_k}{\rightarrow}H^k(B,\Q)$$ vanishes. 

\subsubsection{}
In view of this observation the desirable explanation of the fact
that $\Omega^{k}(\cE)=0$ if  
$k+n\equiv 2(4)$ is that $\ind(\cE)\in K^{-n}(B)$ is in fact of the form
$c_B(\ind_\R(\cE))$ for a real refinement of the index $\ind_\R(\cE)\in
KO^{-n}(B)$.
In fact,  the spinor bundle carries
additional structures which are ``preserved'' upon twisting by real bundles.
Using these structures we can indeed refine the index $\ind(\cE)\in
K^{-n}(B)$ to a class $\ind_\R(\cE)\in KO^{-n}(B)$.

For the purpose of illustration we sketch the construction of
$\ind_\R(\cE)$. Although this is well known, the following exposition
is designed to be a useful reference for the interested reader.

\subsubsection{}
Depending on the class of $n$ modulo $8$ we are going to use quite
different pictures of $KO^{-n}(B)$. We make use of the real Clifford
algebras $C^{p,q}$ associated to $\reals^{p+q}$ with quadratic form
$-x_1^2-\dots- x_p^2 + x_{p+1}^2 +\dots +x^2_{p+q}$.

In one picture an element of $KO^{n}$ is represented as a family of selfadjoint odd Fredholm operators on a graded $C^{n,0}$-module. 
Another representation is as a family of antisymmetric Fredholm operators which anticommute with an action of $C^{0,n-1}$. In this case there is no grading.
Finally an element of $K^1(B)$ is represented by a family of selfadjoint Fredholm operators (and there is again no grading). 
We refer to \cite{atiyahsinger69}  and \cite{karoubi78} for further details.

\subsubsection{$n\equiv 0(8)$\ }

The spinor bundle $\cS(T^v\pi)$ is the complexification of a real
spinor bundle
$\cS_\R(T^v\pi)$. Thus $\cV=\cS(T^v\pi)\otimes W$ is the complexification
of $\cV_\R := \cS_\R(T^v\pi)\otimes W_\R$.
The Dirac operator $D(\cE)$ comes from the Dirac operator $D_\R(\cE)$ on $\cV_\R$. The refined index $\ind_\R(\cE)\in KO^{0}(B)$ is just the index of the family
of real Fredholm operators $D_\R(\cE)^+$.

\subsubsection{$n\equiv 1(8)$\ }

The spinor bundle $\cS(T^v\pi)$ admits a real structure, which anticommutes
with Clifford multiplication. It induces a real structure on $\cV$ which
anticommutes with $D(\cE)$. Let $\cV_\R$ be again the real
$+1$-eigenbundle of the
real structure on $\cV$. The operator $iD(\cE)$ commutes with this real structure and therefore induces an antisymmetric operator $D(\cE)_\R$ on $\cV_\R$.
This family represents $\ind_\R(\cE)\in KO^{-1}(B)$.

\subsubsection{$n\equiv 2(8)$\ }

The spinor bundle $\cS(T^v\pi)$ has a quaternionic structure which is
odd with respect to the $\Z/\Z_2$-grading and commutes with Clifford multiplication. Thus we obtain an induced quaternionic structure $J$ on
$\cV$. We consider $D(\cE)_\R:=JD(\cE)$ as an antisymmetric
real operator on $\cV^+_{\R}$. It anticommutes with the action of
$C^{0,1}$ which is induced by multiplication by $i$.
Therefore the family $D(\cE)_\R$ together with the
$C^{0,1}$-module structure represents $\ind_{\R}(\cE)\in KO^{-2}(B)$.

\subsubsection{$n\equiv 3(8)$\ }

The spinor bundle carries a quaternionic structure which commutes with Clifford multiplication. We get an induced quaternionic structure $J$ on $\cV$ commuting with $D(\cE)$. The antisymmetric operator $D(\cE)_\R:=iD(\cE)$
anticommutes with the action of $C^{0,2}$ generated by $J$ and $iJ$.
Therefore the family $D(\cE)_\R$ on $\cV_\R$ together with the
$C^{0,2}$-module structure represents $\ind_{\R}(\cE)\in KO^{-3}(B)$.

\subsubsection{$n\equiv 4(8)$\ }
 
The spinor bundle $\cS(T^v\pi)$ carries a quaternionic structure which commutes
with the grading and Clifford multiplication.
It induces a quaternionic structure $J$ on $\cV$ which commutes with $D(\cE)$.
 We consider the antisymmetric operator $D(\cE)_\R:=iD(\cE)$
on the bundle $\cV_{\R}$
which anticommutes with the Clifford algebra $C^{0,3}$ generated by
$i,J,iJ$. Therefore the family $D(\cE)_\R$ on $\cV_\R$ together with the
$C^{0,3}$-module structure represents $\ind_{\R}(\cE)\in KO^{-4}(B)$.

\subsubsection{$n\equiv 5(8)$\ }

The spinor bundle $\cS(T^v\pi)$ carries a quaternionic structure
which anticommutes with the Clifford multiplication. It induces a quaternionic structure $J$ on $\cV$ which anticommutes with $D(\cE)$. 
We form the real selfadjoint operator $$D(\cE)_\R:=
\left(\begin{array}{cc}D(\cE)&0\\0&-D(\cE)\end{array}\right)$$ on
$\cV_{\R}\oplus \cV_{\R}$ with its standard odd grading 
$$\left(\begin{array}{cc}0&1\\1&0\end{array}\right)\ .$$
This operator commutes with the Clifford algebra $C^{3,0}$ generated by
$$
\left(\begin{array}{cc}0&J\\-J&0\end{array}\right),\:\left(\begin{array}{cc}0&iJ\\-iJ&0\end{array}\right),\:\left(\begin{array}{cc}1&0\\0&-1\end{array}\right)\ .$$ 
The $C^{3,0}$-equivariant operator $D(\cE)_\R$
represents $\ind_\R(\cE)\in KO^3(B)\cong KO^{-5}(B)$.

\subsubsection{$n\equiv 6(8)$\ }

The spinor bundle $\cS(T^v\pi)$ carries a real structure which
anticommutes with the grading. In induces a real structure $Q$ on
$\cV$ which is odd and commutes with $D(\cE)$. We consider the
selfadjoint operator $D(\cE)_\R:=D(\cE)$ on $\cV_{\R}$. This bundle
is $\Z/2\Z$-graded and admits an action of
$C^{2,0}$ generated by $Q$ and $iQ$ commuting with $D(\cE)_\R$.
The $C^{2,0}$-equivariant operator $D(\cE)_\R$ represents
$\ind(\cE)_\R\in KO^2(B)\cong KO^{-6}(B)$.

\subsubsection{$n\equiv 7(8)$\ }

The spinor bundle $\cS(T^v\pi)$ admits a real structure which commutes with the Clifford multiplication. In induces a real structure $Q$ on $\cV$ which commutes with $D(\cE)$. We consider the real symmetric operator $D(\cE)_\R$ which is obtained by restricting $D(\cE)$ to the $1$-eigenbundle of $Q$.
The operator $D(\cE)_\R$ represents $\ind_\R(\cE)\in KO^{1}(B)\cong KO^{-7}(B)$.

\section{The analytic invariant}

\subsection{Construction of a natural transformation $d^{n}_{B,k-1}$}

\subsubsection{}
Recall that 
complexification of real vector bundles 
induces a natural transformation of multiplicative cohomology theories
$c_B\colon KO^n(B)\rightarrow K^n(B)$. The real index
$\ind_\R(\cE)\in KO^{-n}(B)$ is a refinement
of $\ind(\cE)\in K^{-n}(B)$ in the sense that 
$$c_B(\ind_\R(\cE))=\ind(\cE)\ .$$

\subsubsection{}
For $k\ge 0$ and $n\in\Z$
we define the group $U_k^n(B)$ by the following exact sequence
$$0\rightarrow U^n_k(B)\rightarrow KO^n(B)\stackrel{q_B}{\rightarrow} K^n(B)/K_k^n(B)\ ,$$
where $q_B$ is the composition of $c_B$ with  the projection onto the quotient. We
also define $U^n_{\infty}(B)$ by the exact sequence
\begin{equation*}
  0 \to U_\infty^n(B) \to KO^n(B) \to K^n(B).
\end{equation*}
Since $q_B$ is a natural transformation the association $B\mapsto
U^n_k(B)$ extends to a functor with values in abelian groups.

\subsubsection{}
Assume that $k-n\equiv 2(4)$.
\begin{definition}\label{def:of natural transformation}
We  define the natural transformation
  \begin{equation*}
d_{B,k-1}^n\colon U^n_k(B)\rightarrow H^{k-1}(B,\R/\Z)
  \end{equation*}
  by the requirement that $$d_{B,k-1}^n(x):=\hat c_{k-1}(\cE)\ ,$$
where $\cE$ is any geometric family of dimension $8l-n$, $l\in\Z$, such that $x=\ind_\R(\cE)$.
 \end{definition}

%\begin{remark}
%  This natural transformation, defined analytically with the help of
%  local index theory, is the main object of study of the rest of this
%  paper. In particular, we will identify it in topological terms, and
%  will prove vanishing and non-vanishing results (in dependence on $k$
%  and $n$).
%\end{remark}

\subsubsection{}
We must check that  the definition of $d_{B,k-1}^n $ makes sense. 
\begin{prop}\label{plam}
If $k-n\equiv 2(4)$, then there is a unique homomorphism
$$d_{B,k-1}^n\colon U^n_k(B)\rightarrow H^{k-1}(B,\R/\Z)$$ such that
$d_{B,k-1}^n(\ind_\R(\cE))=\hat c_{k}(\cE)$ for any 
geometric family over $B$ of dimension $8l-n$, $l\in\Z$, with $\ind_\R(\cE)\in U^n_k(B)$.
This homomorphism is
natural with respect to continuous maps.
\end{prop}
\begin{proof} The essential parts of the proof are given in Lemma \ref{l1}
and Lemma \ref{l2} below.
It immediately follows from these Lemmas that for given $B$ there exist a unique
map $d^n_{B,k-1}$ satisfying the requirement.
Additivity and naturality with respect to smooth maps of $d^n_{B,k-1}$
follows from naturality and additivity of the class  $\hat
c_{k}(\cE)$.
But then naturality extends to continuous maps since
 $U^n_k(\cdot)$ as well as
$H^{k-1}(\cdot,\reals/\integers)$ are weak homotopy functors. \end{proof}

\subsubsection{}
\begin{lem}\label{l1}
If $\ind_\R(\cE)=0$, then $\hat c_{k}(\cE)=0$.
\end{lem}
\begin{proof} 
Assume that $\ind_\R(\cE)=0$. In this case we can find a smoothing
perturbation of the real operator $D(\cE)_\R$ which is invertible.
We call this perturbation a real taming.
By complexification a real taming
induces a taming $\cE_t$ which is compatible with the additional
symmetries determining the real structure.

These additional symmetries imply that the Chern form of the
Bismut super connection associated to $D(\cE)$ and its tamed
perturbation vanishes.  Since the $\eta$-form is derived from this Chern form we
conclude that $\eta(\cE_t)=0$ if the taming is induced from a
real taming.  
The assertion of the Lemma now follows from the description of $\hat
c_k(\cE)$ in terms of the $\eta$-form (see Subsection \ref{co2}). \end{proof} 

\subsubsection{}
\begin{lem}\label{l2}
If $x\in KO^n(B)$, then there exists a geometric family $\cE$ as above
such that
$\ind_\R(\cE)=x$.
\end{lem}
\begin{proof}   
By the periodicity of $KO^*(B)$ we can assume that $n<0$. 
By definition, $${KO}^{n}(B) =
\widetilde{KO}^{n}(B_+)=\widetilde{KO}^0(\Sigma^n
B_+)\ ,$$ where $B_+$ is obtained from $B$ by adjoining 
an additional base point, and $ \widetilde{KO}^{n}(B_+)$ denotes the
reduced $KO$-theory.

Let now $x\in {KO}^{n}(B)$ correspond to
$\tilde x\in \widetilde{KO}^0(\Sigma^n
B_+)$.
Let  $p\colon S^n\times B_+\rightarrow \Sigma^n B_+$ be the natural
projection and $W_\R=W_\R^+\oplus W_\R^-$ be the real $\Z/2\Z$-graded vector bundle over $S^n\times
B_+$ representing $p^*\tilde x$.

We form two geometric families $\cE^\pm$ with
underlying bundle
$S^n\times B_+\rightarrow B_+$ with its standard fibrewise orientation
and spin structure,  and with the real twisting bundle $W_\R^{\pm}$.
Then we have (by Bott periodicity or the index theorem)
$x=\ind_\R(\cE^+\cup_{B_+}\cE^-)_{|_B}$.
\end{proof} 

\subsection{Some properties of $d^n_{B,k-1}$}

%The main goal of the present paper is to give a topological interpretation of the transformations $d^n_{B,k-1}$. First we collect some simple consequences
%of the analytic construction.

\subsubsection{}\label{forhand}
We approach the study of the natural transformation
$d^n_{B,k-1}$ from two sides. First, in view of its definition through
the analysis of families of Dirac operators we use mainly analytical arguments
in order to show some simple properties of this transformation.
This is the subject to the present section. 

A finer study in Section \ref{univsec} leading to a complete understanding of the transformation
uses methods from topology and the observation that a natural
transformation comes from an universal one between suitable classifying
spaces. It should be noted that most results of this section, the
important exception being Corollary \ref{jh},  will also
follow from the topological description, and will not be needed to
derive this description.

\subsubsection{} 
There is a natural transformation $$r_B:K^*(B)\rightarrow KO^*(B)\ .$$
It is determined by the special case $r_B:K^0(B)\rightarrow KO^0(B)$ 
which associates to a class represented by a complex vector bundle the 
class represented by the underlying real vector bundle.
It is easy to see that 
$$r_B\circ c_B=2$$
(multiplication by $2$).
\subsubsection{}

\begin{prop}\label{zmal0}
We have $2d^n_{B,k-1}=0$.
\end{prop}
\begin{proof}
Fix $x\in U^n_k(B)$. The homology class $d^n_{B,k-1}(x)\in
H^{k-1}(B,\R/\Z)$ is determined by its values on all smooth cycles $Z\in
Z^{k-1}$ on $B$. 

Given a $(k-1)$-cycle $Z$ there exists a manifold $X$ which is 
homotopy
equivalent to a $(k-1)$-dimensional
CW-complex, a smooth map $f:X\rightarrow B$, and a $(k-1)$-cycle
$Z^\prime$ in $X$ such that $f_*Z^\prime=Z$ (compare \ref{thik}).
By the naturality of $d^{n}_{\cdot,k-1}$ we have 
$$2 d^n_{B,k-1}(x)(Z)=2
f^*d^n_{B,k-1}(x)(Z^\prime)=d^n_{X,k-1}(2f^*x)(Z^\prime)\ .$$
It thus suffices to show that $2f^*x=0$.

Since $f^*x\in U_k^n(X)$ and $X$ is (up to homotopy equivalence) $(k-1)$-dimensional
we have  $c_X(f^*x)=0$. This implies $0=r_X\circ c_X(f^*x)=2f^*x$.
\end{proof}

%Let $\Sigma$ be a arbitrary manifold homotopy
%equivalent to a $(k-1)$-dimensional
%CW-complexes and $f\colon \Sigma\to B$ an arbitrary map. 
%%Let $f\colon \Sigma\rightarrow B$ and $a\in H_{k-1}(\Sigma;\Z)$, giving rise
%%to a class $f_*(a)\in H_{k-1}(B,\Z)$. 
%As explained in Subsection
%\ref{co2} and because of the naturality of $d^n_{\cdot,k}$
%it suffices to show that $2d^n_{\Sigma,k-1}(f^*x)=0$.
%Note that $f^*x\in U^n_{\infty}(\Sigma)$ for dimension reasons, which
%implies that $c_\Sigma(f^*x)=0$.
%Now there is a natural transformation
%$r_\Sigma\colon K^n(\Sigma)\rightarrow KO^n(\Sigma)$ which is induced by
%forgetting the complex  structure. It satisfies $r_\Sigma\circ
%c_\Sigma=2$ (multiplication with $2$).
%Thus $2f^*x=0$. Therefore,
%$2d^n_{\Sigma,k-1}(f^*x)=0$. \end{proof} 

\subsubsection{}
We have defined $d^n_{B,k-1}$ in terms of eta forms of families of
Dirac operators. In the following we show that it suffices to
understand eta forms for zero-dimensional families. Note that the
latter are essentially objects of linear algebra.

\subsubsection{}

If $x\in U^{-n}_\infty(B)$, then there exist a real $\Z/2\Z$-graded
vector bundle 
over $\Sigma^nB_+$ which represents the class $\tilde x\in
\widetilde{KO}^0(\Sigma^nB_+)$ corresponding to $x$ under the identification
$KO^{-n}(B) \cong \widetilde{KO}^0(\Sigma^n B_+)$.
Let $W_\R=W^+_{\R}\oplus W^-_{\R}$ be the pull-back of this bundle to $S^n\times B$
under the natural map $S^n\times B\rightarrow \Sigma^nB_+$.

The bundle $\pi\colon S^n\times
B\rightarrow B$ has a natural fibrewise orientation and spin
structure. The round metric of $S^n$ induces a fibrewise metric $g^{T^v\pi}$.
The canonical decomposition $T(S^n\times B)=\pr_{S^n}^*TS^n\oplus \pr_B^*TB$
yields the horizontal distribution $\pr_B^*TB$.
After choosing geometric  
bundles 
$\bW_\R^\pm=(W^\pm_\R,\nabla^{W^\pm_\R},h^{W^\pm_\R})$ we
obtain the geometric family $\cE^\pm$ and $\cE:=\cE^+\cup_{B}
(\cE^{-})^{ op}$ over $B$, with underlying bundle $S^n\times B$, such that $\ind_\R(\cE)=x$.

\subsubsection{} 
Since $c_B(x)=0$  by assumption we know that the complexification $W$ of $W_\R$ represents the trivial element
$0=[W]\in K^0(S^n\times B)$. Thus (possibly after adding a trivial
bundle of formal dimension zero) the bundle $\bW$ admits an odd unitary
selfadjoint (not necessary parallel) automorphism $U$.

\subsubsection{}
The bundles $\bW_\R^\pm$ give rise to geometric families $\cF^\pm$ over
$S^n\times B$ where the underlying zero-dimensional fiber bundle is $\id\colon S^n\times
B\rightarrow S^n\times B$, and the twisting bundle is
$\bW_\R^\pm$.
We let $\cF:=\cF^+\cup_{S^n\times B}(\cF^-)^{op}$.
Then we have $\ind_\R(\cF)=[W_\R]\in KO^0(S^n\times B)$
and $\ind(\cF)=[W]=0\in K^0(S^n\times B)$.
The automorphism $U$ gives a taming $\cF_t$ of $\cF$.
In particular, we have a well-defined form $\eta(\cF_t)\in
\cA(S^n\times B)$ such that
$d\eta(\cF_t)=\Omega(\cF)=\ch(\nabla^{W})$.

\subsubsection{}
For $r>>0$ the operator $rU$ can be considered as a sort of taming of the family $\cE$. It is not a taming
in the strong sense since $U$ is not smoothing along the
fibers of $\pi$. Rather it is a local taming 
in the sense of \cite{bunke020}. Local index theory works for local
tamings as well. 
We let $\cE_t(r)$ be the geometric family $\cE$ tamed with $rU$
and let $\eta(\cE_t(r))$ be the associated $\eta$-form.
\subsubsection{}
\begin{prop}\label{fe}
We have
$$\lim_{r\to\infty} \eta(\cE_t(r))=\int_{S^n\times B/B} \eta(\cF)$$
\end{prop}
\begin{proof}
This assertion is proved using the adiabatic limit techniques developed e.g. in
\cite{bunkema02}. The general method gives
$$\lim_{r\to\infty} \eta(\cE_t(r))=\int_{S^n\times B/B} \hA(\nabla^{T^v\pi}) \eta(\cF_t)\ .$$
The result now follows from $
\hA(\nabla^{T^v\pi})=\pr_{S^n}^*\hA(\nabla^{TS^n})=1$, since
$\hA(\nabla^{TS^n})=1$ for round metric (which is also locally
conformally flat). \end{proof}

\subsubsection{}
A closed form $\omega\in \cA^{k-1}(B)$ represents a class
$[\omega]\in H^{k-1}(B,\R)$. By $[\omega]_{\R/\Z}\in H^{k-1}(B,\R/\Z)$
we denote its natural image.
Let $m$ be determined by $2m=k$ or $2m=k+1$.
Then definition \eqref{eq:compute_char_class} together with
Proposition \ref{fe} implies
 the following corollary.
\begin{corollary}\label{jh}
With the notation above
$$d^{-n}_{B,k-1}(x)=[(-1)^{m-1}(m-1)! \int_{S^n\times B/B} \eta^{k+n-1}(\cF)]_{\R/\Z}\ .$$\end{corollary}
%where $[\omega]_{\R/\Z}\in H^{k-1}(B,\R/Z)$
%is the image class represented by the closed form $\omega$
%under the map $H^*_{dR}(B)\cong H^*(B,\R)\rightarrow H^*(B,\R/\Z)$ and
%$m$ is determined by $2m=k$ or $2m=k+1$.
In fact, in view of Proposition \ref{zmal0} we could also omit the
sign $(-1)^{m-1}$.
\subsubsection{} 
We consider the sequence
$$0\rightarrow \Z\rightarrow \R\rightarrow \R/\Z\rightarrow 0$$
and let $\beta^{\prime\prime}\colon H^{k-1}(B,\R/\Z)\rightarrow
H^k(B,\Z)$ be the associated Bockstein operator. Recall from \ref{firstbock}, that
$\beta^{\prime\prime}$ coincides with the composition 
$$H^{k-1}(B,\R/\Z)\rightarrow H_{Del}^{k}(B)\stackrel{v}{\rightarrow} H^{k}(B,\Z)\ .$$
\begin{prop}
For $x\in U^0_{4k+2}(B)$ we have $$(\beta^{\prime\prime}\circ
d^0_{4k+1})(x)=(c_{4k+2}\circ c_B)(x)\ .$$\end{prop}
\begin{proof}
We write $x=\ind_\R(\cE)$ for a suitable geometric family $\cE$ over
$B$. Then
we have the chain of equalities
$$(\beta^{\prime\prime}\circ d^0_{4k+1})(x)=(v\circ \hat
c_{4k+2})(\cE)=c_{4k+2}(\ind(\cE))=(c_{4k+2}\circ c_B)(x)\ .$$ 
\end{proof}

\subsection{Nontrivial examples}\label{trio}

In Section \ref{sec:transgression} we will give a complete description
of the universal classes
$d^n_{\cdot,k}$, which also decides when exactly these classes can be
non-trivial. In this section, we want to construct explicit and easy
non-trivial examples over low dimensional spheres as base manifolds.

\subsubsection{}\label{surd}
Let $MSpin^*$ be the spin bordism cohomology theory and
$\alpha:MSpin^*\rightarrow KO^*$ be the $\alpha$-genus
introduced by Hitchin \cite{hitchin74}.
Note that $\alpha:MSpin^*(*)\rightarrow KO^*(*)$
is surjective. If $E$ is a closed spin manifold, then we write
$\alpha(E)$ for the result of $\alpha$ applied to the class
$[E]\in MSpin^*(*)$
represented by $E$

The most important common feature of the following examples is
that they all come with a trivial twisting bundle.
In other words, the respective geometric family $\cE$ represents an
element  $[E,\pi]\in MSpin^{-n}(S^i)$ (with $n\in\nat_0$  and $i\in\{0,1,2\}$
depending on the case), such that $\ind_\R(\cE)=\alpha([E,\pi])$.

\subsubsection{}{[ $n\equiv1(8)$, $k=1$]}\label{lmk}

\newcommand{\MSpin}{\mathbf{MSpin}}
 Let $E$ be a closed  spin manifold of dimension
$n\equiv 1(8)$  with $\alpha(E)=1\in KO^{-1}(*)\cong \Z/2\Z$.
Such manifolds exist by the Remark \ref{surd}. 

We choose a Riemannian metric $g^{TE}$.
Then we consider $E$ as a geometric family $\cE$ over the point $*$
with the trivial twisting bundle $W_\R=E\times\R$.
We claim that $\hat c_1(\cE)\not=0$. 

Note that
$H^1_{Del}(*)=\R/\Z$ and $\hat c_1(\cE)([*])=[ \eta^0(\cE_t)]_{\R/\Z}$,
where $\cE_t$ is any taming. The degree $0$ part $\eta^0(\cE_t)$ is defined
even for a pre-tamed manifold in the sense of \cite{bunke020}, and if
the pre-taming is trivial, then it is
one half of the usual Atiyah-Patodi-Singer eta invariant \cite{atiyahpatodisinger75}. 
 Now we have $\eta^0(\cE)=0$ since $D(\cE)$ anticommutes with
the real structure and thus has symmetric spectrum. 
Since by spectral flow (taming essentially means that the underlying operators are invertible)
\begin{equation*}
[\eta^0(\cE_t)-\eta^0(\cE)]_{\R/\Z}=[\frac{1}{2} \dim\ker
D(\cE)]_{\R/\Z}
\end{equation*}
we see that
$$[\eta^0(\cE_t)]_{\R/\Z} = [\frac12 \dim\ker D(\cE)]_{\R/\Z}\ .$$
The condition $\alpha(E)\not=0$ says that $\dim\ker D(\cE) \equiv 1(2)$.
Therefore,
$d_{*,0}^1(\ind_\R(\cE))=[\eta^0(\cE_t)]_{\R/\Z}=[\frac12]_{\R/\Z}\not=0\in
H^0(*,\R/\Z)$.

\subsubsection{}{[$n\equiv 0(8)$, $k=2$]}\label{kml1}

We consider a family $E\rightarrow S^1$  of closed spin manifolds with fiber
dimension
$n\equiv 0(8)$ and 
$\alpha(E)\not=0$. Indeed for  any given spin manifold $M$ of dimension $n\equiv 0(8)$
such a bundle with fiber $M$ exists
by \cite{hitchin74}. We will in addition assume that $\alpha(M)\stackrel{\text{Def}}{=}\hat{A}(M)=0$.

We equip $E$ with geometric structures and consider the trivial
twisting bundle $E\times W_\R$. Let $\cE$ denote the corresponding
geometric family over $S^1$. Since $\alpha(M)=0$ we have
$\ind(\cE)\in K_2^0(S^1)=0$ and thus $\ind_\R(\cE)\in U_2^0(S^1)$.

We claim that $\hat c_2(\cE)\not=0$.
We consider a taming $\cE_t$. For $a>0$ the standard metric $g^{TS^1}$ of $S^1$ and the
horizontal distribution of $\cE$ induces a Riemannian metric
$g_a^{TE}=g^{T^v\pi}\oplus a \pi^*g^{TS^1}$ on the
total space $E$. Since $E$ is spin we can consider the total Dirac operator
$D(a)$ on $E$  and its perturbation $D_t(a)$ which is induced by the taming.
In the adiabatic limit $a\to 0$ the operator $D_t(a)$ becomes
invertible.
In other words, for small $a$ the perturbation $D_t(a)$ is induced by
a local taming. 
As in \cite{bismutcheeger89}, we have
$$[\lim_{a\to 0} \eta^0(D_t(a))]_{\R/\Z}= [\int_{S^1}\eta^1(\cE_t)]_{\R/\Z} =d^n_{S^1,1}(\ind_\R(\cE))([S^1])\in\R/\Z\ .$$ 
For sufficiently small $a$ the class $[\eta^0(D_t)(a)]_{\R/\Z}$ is independent
of the adiabatic parameter. As in \ref{lmk}, since $\alpha(E)\not=0$
and $\eta^0(D(a))=0$ we have $[\eta^0(D_t(a))]_{\R/\Z}=
[\frac12]_{\R/\Z}$. Thus $d^n_{S^1,1}\not=0\in H^1(S^1,\R/\Z)$.
 
This result can be interpreted as follows : The holonomy of the
determinant line bundle of $\cE$ is $-1\in U(1)$.

\subsubsection{}{[$n=7(8)$, $k=3$]}\label{sec:ex_n_equal_7}

We  consider a family  $E\rightarrow S^2$ of closed spin manifolds
with fiber dimension $n\equiv 7(8)$ and  $\alpha(E)\not=0$.
Such a family exists \cite{hitchin74} for any given closed spin manifold
$M$ of dimension $n\equiv 7(8)$. We choose geometric structures and
consider the trivial twisting bundle $E\times\R$.
In this way we obtain a geometric family $\cE$.
Since $K^1(S^2)=0$ we have $\ind(\cE)=0$ and therefore
$\ind_\R(\cE)\in U_3^7(S^2)$.
We claim, that
$d^7_{S^2,2}(\ind_\R(\cE))\not=0\in H^2(S^2,\R/\Z)$. 

We proceed as in \ref{kml1}.
We consider a taming $\cE_t$. It induces a perturbation
$D_t(a)$ of the total operator $D(a)$ on $E$.
We have again
$$[\lim_{a\to 0} \eta^0(D_t(a))]_{\R/\Z}=[ \int_{S^2}\eta^2(\cE_t)]_{\R/\Z}= d^7_{S^2,2}(\ind_\R(\cE))([S^2])  \ .$$ 
Again, for sufficiently small $a$ the class  $[\eta^0(D_t(a))]_{\R/\Z}$ is independent of
the adiabatic parameter. Since $\alpha(E)\not=0$ and $\eta^0(D(a))=0$
we have $[\eta^0(D_t(a))]_{\R/\Z} = [\frac12]_{\R/\Z}$. 
This implies the claim.

%\Kommentar{This should be explained in somewhat more detail!}

\section{Topological universal classes}\label{univsec}

\subsection{Transgression of the Chern classes}

\subsubsection{}
As proposed in \ref{forhand} we shall understand $d^0_{B,4k+1}$
through its universal example. In the present section we start with
the definition of this universal class. We will obtain an expression
of this class in terms of familiar characteristic classes of real
vector bundles. In Theorem \ref{top1} we show that the transformation
$d^0_{B,4k+1}$ is indeed induced by the corresponding universal example.
In Theorem \ref{theo:top2} we will then identify
$d^{-n}_{B,4k+1-n}$ in topological terms for all $n$.
Then we will establish some vanishing results and more details about the
topological side.

%We will show later that our natural transformation of Definition
%\ref{def:of natural transformation} is defined by pulling back a
%universal cohomology class from a space which plays the role of a
%classifying space for elements in $U^n_k$, and that this universal
%cohomology class, in turn, can be expressed in terms of familiar
%cohomology classes of vector bundles. In this section, we first define
%and study this universal cohomology class, and then proof that the
%corresponding natural transformation is indeed the analytically
%defined natural transformation of Definition \ref{def:of natural
%  transformation}.

\subsubsection{} 

In the present section all spaces have distinguished base points and all maps are base point preserving.
Let $O$ and $U$ be the direct limits of $O(n)$ and $U(n)$ induced by
the embedding into the left upper corner. The embeddings
$\R^n\hookrightarrow \C^n$, $n\in\nat$, induce embeddings
$O(n)\hookrightarrow U(n)$ and $c\colon O\hookrightarrow U$.
The map $c$ induces the complexification transformation $c_B$.
(see \ref{cbintro}).

\subsubsection{}
Let $EU\rightarrow BU$ be a universal bundle for $U$.
We can consider $BO:=EU/O$ and obtain a bundle
\begin{equation}\label{eq:U_by_O_fibration}
U/O\xrightarrow{i} BO\xrightarrow{p} BU
\end{equation}
with fiber $U/O$ over the base point of $BU$, compare
\ref{item:projection BO to BU} and \ref{item:U_modO_toBU} on page 
\pageref{item:U_modO_toBU}.

In the following, we use the transgression homomorphism for
cohomology associated to this fibration \eqref{eq:U_by_O_fibration}.
 For the convenience of the reader, we have
collected the main definitions and properties of transgression in
general (with proofs) in
Appendix \ref{sec:appendix_transgresseion}, and of transgression and
cohomology of the spaces in the fibration $U/O\to BO\to BU$ in Lemma \ref{lem:transgressive classes}
in Appendix \ref{sec:cohomology-bo-bu}.

\subsubsection{}
\begin{ddd}
%\begin{enumerate}
%\item For $x\in H^{4k+2}(BU,\Q)$ we define its transgression
%$T(x):=T_h(x)\in H^{4k+2}(U/O,\Q)$
%(for any choice of $h$).
%\item 
We define the \emph{universal transgressed Chern classes}
  $$d_{4k+1}:=T(c^\Q_{4k+2})\in H^{4k+1}(U/O,\Q)\ ,$$
where $c^\Q_{4k+2}$ is the image of the universal Chern class
$c_{4k+2}$ under the natural map
\begin{equation*}
H^{4k+2}(BU,\Z)\rightarrow H^{4k+2}(BU,\Q).
\end{equation*}
%\end{enumerate}
\end{ddd}

%\begin{lem}\label{deco}
%If $x\in H^{4k+2}(BU,\Q)$  can be decomposed as $x=y\cup z$
%with $y\in H^{4l+2}(BU,\Q)$ and $z\in H^{4(k-l)}(BU,\Q)$ and $l<k$, then we have $T(x)=0$.
%\end{lem}
%\begin{proof}
%Let $P:K(\Q,4l+2)\wedge K(\Q,4(k-l))\rightarrow K(\Q,4k+2)$ be the map representing the
%cup product. Then $x=P\circ(y\times z):BU\rightarrow K(\Q,4k+2)$.
%Let $h:CBO\rightarrow K(\Q,4l+2)$ be the extension (\ref{hh1}) as constructed above for
%$y$. Then $u^*(h)\times (\pr_{BO}^* \circ p^*)(z):I\times BO\rightarrow K(\Q,4l+2)\wedge K(\Q,4(k-l))$ factors over  the quotient  $u:I\times BO\rightarrow CBO$.
%Let $\tilde h$ denote this factorization.
%Then $P\circ \tilde h:CBO\rightarrow  K(\Q,4k+2)$ is the homotopy (\ref{hh1})
%for $x$. But $\tilde h\circ C(i)$ is  constant. Now the assertion follows
%from the construction of $T(x)$. \end{proof}

\subsubsection{}\label{lsdeef}
We now consider the following commutative diagram:
\begin{equation}
\begin{CD}
0 @>>>\Z @>{2\times}>>\Z @>{q}>>  \Z/2\Z @>>>0\\
&& @VV{=}V @VV{\frac{1}{2}\times}V  @VV{l}V\\   
0 @>>> \Z @>>>  \Q @>>> \Q/\Z @>>> 0
\end{CD}\ .\label{eq:def_of l}
\end{equation}
 
\begin{ddd}
We define $\bar d_{4k+1}\in H^{4k+1}(U/O,\Q/\Z)$ as the image of $d_{4k+1}$
under the natural map $H^{4k+1}(U/O,\Q)\rightarrow H^{4k+1}(U/O,\Q/\Z)$.
\end{ddd}

\subsubsection{}
The cohomology ring of $BO$ with coefficients in $\Z/2\Z$ is a
polynomial ring $$\Z/2\Z[w_1,w_2,\dots]\ ,$$ where $w_i\in
H^{i}(BO,\Z/2\Z)$ are the universal Stiefel-Whitney classes.
It is well known (see \cite{milnorstasheff68})
that 
\begin{equation}\label{cb}\beta(w_{2k}\cup
  w_{2k+1})=p^*c_{4k+2},\end{equation}
where $\beta$ is the cohomological Bockstein operator associated to the exact
sequence of coefficients in the first row of \eqref{eq:def_of
  l}. In particular, $2p^*c_{4k+2}=0 \in H^{4k+2}(BO;\integers)$, as
also stated in Appendix \ref{sec:cohomology-bo-bu}.
 
\subsubsection{}
Let $l_*:H^{4k+1}(BO,\Z/2\Z)\rightarrow H^{4k+1}(BO,\Q/\Z)$ be induced by
$l$ of \eqref{eq:def_of l}. Note that by Lemma \ref{lem:l factors
  through Bockstein} the map $l_*$ factors over the image of $\beta$.
 \begin{definition}\label{def:of tilde d}
 We define 
  \begin{equation*}
\tilde d_{4k+1} := l_*u \in H^{4k+1}(BO,\Q/\Z),
\end{equation*}
where  $u\in H^{4k+1}(BO,\Z/2\Z)$ is such that $\beta(u)=p^*
  c_{4k+2}$.
\end{definition}
  Such an $u$ exists by (\ref{cb}) and $l_*u$ is independent of the
  choice since we have fixed $\beta(u)$.
  
\subsubsection{}

Let $i\colon U/O\hookrightarrow BO$ be the inclusion.
\begin{lem}
We have $i^* \tilde d_{4k+1} = \bar d_{4k+1}$.  
Moreover, $2\tilde d_{4k+1}=0$ and $2\bar d_{4k+1}=0$.
\end{lem}
\begin{proof}
The first assertion is a special case of Proposition 
\ref{prop:transgression and
  Bockstein}, where we use $2p^*(c_{4k+2})=0$. 
Note that the homomorphism $l$
is given by division by $2$.

The second assertion follows from the fact that $2\tilde d_{4k+1}= 2l_*
u = l_*(2u)=0$ for $u$ of Definition \ref{def:of tilde d}, and
$2\overline{d}_{4k+1}= i^*(2\tilde d_{4k+1})=0$. 
\end{proof} 

\subsubsection{}
By \eqref{cb} we have the following corollary. 
\begin{kor}\label{spi1}
We have $\tilde d_{4k+1}=l_* (w_{2k}\cup w_{2k+1})$.
\end{kor}

\subsection{A topological description of $d^0_{B,4k+1}$}

\subsubsection{}
In this subsection we pretend that $p\colon BO\rightarrow BU$ is a smooth
fiber bundle. To be precise, we should replace this bundle by a
$N$-equivalent finite-dimensional smooth bundle for $N$ sufficiently
large.
 
Let $W_\R^+\rightarrow BO$ be the universal bundle.
Then $p\colon BO\rightarrow BU$ classifies its  complexification, i.e. if
$W^+\rightarrow BU$ is the universal bundle over $BU$, then
we have an isomorphism $W_\R^+\otimes_\R\C\cong p^*W^+$ which induces a real structure $Q$ (complex conjugation) on $p^*W^+$.
We can assume that $W_\R^+$ comes with a metric $h^{W_\R^+}$, and we
choose a connection $\nabla^{W_\R^+}$.
We set $\bW_\R^+:=(W_\R^+,\nabla^{W_\R^+},h^{W_\R^+})$
and let $\bW^+$ be its complexification.

\subsubsection{}
We now consider the $\Z/2\Z$-graded bundle $\bW:=\bW^+\oplus  \bW^-$ with
$\bW^-:=\overline{ \bW}^+$. 
It admits an odd unitary selfadjoint (not necessary parallel) automorphism
$$R:=\left(\begin{array}{cc}0& Q\\ Q & 0\end{array}\right)\ .$$
We form the geometric families $\cG^\pm$ on $BO$  with underlying fiber bundle
$\id\colon B\rightarrow B$ and twisting bundles $p^*\bW^{\pm}$. Then the family
$\cG=\cG^+\cup_B \cG^-$ admits a taming $\cG_t$ induced by $R$.
The associated $\eta$-form satisfies
$d\eta(\cG_t)=\Omega(\cG)=\ch(\nabla^{p^*W})$.
By construction we have
$\ch_{4k+2}(\nabla^{p^*W})=2p^*\ch_{4k+2}(\nabla^{W^+})$ and
$\ch_{4k}(\nabla^{p^*W})=0$ (compare \ref{req1}).

\subsubsection{}
Let $i\colon U/O\rightarrow BO$ be the inclusion of the fiber.
Then
\begin{equation*}
d i^*\eta^{4k+1}(\cG_t)=2(i^*\circ p^*)\ch_{4k+2}(\nabla^{W^+})=0
\end{equation*}
since $p\circ i$ is a constant map. Thus we can consider the class
\begin{equation*}
[i^*\eta^{4k+1}(\cG_t)]\in H^{4k+1}(U/O,\R).
\end{equation*}

\subsubsection{}
Let $d^\R_{4k+1}\in H^{4k+1}(U/O,\R)$ be the 
image of $d_{4k+1}\in H^{4k+1}(U/O,\Q)$ under the natural map
induced by the map of coefficients $r\colon \Q\into \R$.
\begin{lem}\label{coonex}
We have 
$$[i^*\eta^{4k+1}(\cG_t)]=  \frac{2d_{4k+1}^\R}{(2k)!}  \ .$$
\end{lem}
\begin{proof}
The proof follows from the fact that in the smooth situation there is
the alternative description of the transgression
$T^\R\colon H^{4k+2}(BU,\R)\rightarrow H^{4k+1}(U/O,\R)$ given in
Proposition \ref{prop:explicit_transgression}. Let $x\in
H^{4k+2}(BU,\R)$ be represented by a closed form
$X\in \cA^{4k+2}(BU)$. Then there is some form $Y\in \cA^{4k+1}(BO)$
such that $dY=p^*X$. The class
$T^\R(x)\in H_{dR}^{4k+1}(U/O,\R)$ is then represented by $i^*Y\in \cA^{4k+1}(U/O)$. 

In our case
\begin{equation*}
x=2\ch^\R_{4k+2}([W^+]),\quad X=2\ch_{4k+2}(\nabla^{W^+}), \quad
Y=\eta^{4k+1}(\cG_t)
\end{equation*}
 so that
 \begin{equation*}
[i^*\eta^{4k+1}(\cG_t)]=T^\R(2\ch^\R_{4k+2}([W^+])).
\end{equation*}
Note that 
$$\ch_{4k+2}([W^+])=\frac{1}{(2k)!} c^\Q_{4k+2}+\mbox{decomposable
  classes},$$
where for reasons of degree each decomposable summand contains at
least one factor $c_{4l+2}$ which is transgressive by the table on
page \ref{table:maps-between-loop their effect in homology}. 
Hence by Proposition \ref{prop:transgression_of_products}
$$T(2\ch_{4k+2}([W^+]))=\frac{2}{(2k)!} T(c_{4k+2}) = \frac{2 d_{4k+1}}{(2k)!} \ .$$
This implies the assertion since
$T^\R\circ r_* = r_*\circ T$ by Lemma \ref{lem:transgression is
  compatible with natural transformations}.    \end{proof} 

\subsubsection{}
We now consider a manifold $B$ and $x\in \widetilde{KO}^0(B)$.
Let $X\colon B\rightarrow BO$
be a classifying map for $x$.
We assume that $x\in U^0_\infty(B)$. Then we can assume that
$X$ factors through the inclusion $i\colon U/O\rightarrow BO$, i.e.
without loss of generality we can assume that $X\colon B\rightarrow U/O$.

We  define $\bar d^\R_{4k+1}\in H^{4k+1}(B,\R/\Z)$
as the image of $d^\R_{4k+1}$ under the map of coefficients
$\R\rightarrow \R/\Z$, or equivalently, as the image of
$\bar d_{4k+1}$ under the map of coefficients
$\Q/\Z\rightarrow \R/\Z$. 
We now come to the main result of this subsection.
\begin{theorem}\label{top1}
$d^0_{B,4k+1}(x)= X^*\bar d^\R_{4k+1}$.
\end{theorem}
\begin{proof}
Let $x$ be of the form $[V_\R^+]-[V_\R^-]$, where $V_\R^-:=\R^N\times B$
is trivial. Then we have an isomorphism $V_\R^+\cong X^*W_\R^+$. 
The metric $h^{W^+_\R}$ and the connection $\nabla^{W_\R^+}$ induce a
metric $h^{V^+_\R}$ and  a connection $\nabla^{V_\R^+}$ on $V_\R^+$. In this way
we obtain a geometric bundle
$\bV^+_\R=(V_\R^+,\nabla^{V_\R^+},h^{V^+_\R})$. Furthermore,
we equip $V^-_\R$ with the canonical geometry and get $\bV^-_\R$.

Set $\bV^\pm:=\bV_\R^\pm\otimes \C$ and consider the $\Z/2\Z$-graded
bundle $\bV:=\bV^+\oplus \bV^-$. Since $[V]=0$ in $K^0(B)$ we can choose a 
unitary odd selfadjoint (not necessary parallel) automorphism $U$ of $V$. 

We form the geometric families $\cH^\pm$ over $B$ with underlying bundle
$\id\colon B\rightarrow B$ and twisting bundle $\bV^\pm$. Furthermore we define 
$\cH:=\cH^+\cup_B(\cH^-)^{op}$. Then we have
$\ind_{\R}(\cH)=x$
and 
$d^0_{4k+1,B}(x)=\hat c_{4k+2}(\cH)$.
The isomorphism $U$ induces a taming $\cH_t$.
By Corollary \ref{jh} we thus have 
$$\hat c_{4k+2}(\cH) =  [(2k)!\eta^{4k+1}(\cH_t)]_{\R/\Z}\in
H^{4k+1}(B,\R/\Z)\ .$$

We now consider the bundle $\tilde \bV:=\bV^+\oplus \bar \bV^-\oplus \bar \bV^+\oplus \bV^-$ with the $\Z/2\Z$-grading $\diag(1,1,-1,-1)$ and the two odd automorphisms
$$\tilde U:=\left(\begin{array}{cccc}0&0&0&U^-\\0&0&\bar U^+&0\\
0&\bar U^-&0&0\\U^+&0&0&0\end{array}\right),\quad \tilde R:=\left(\begin{array}{cccc}0&0&R^+&0\\0&0&0&R^-\\
R^+&0&0&0\\0 &R^-&0&0\end{array}\right)\ ,
$$
where $R^\pm$ is the $\complexs$-linear isomorphism between
$V^\pm_\reals\tensor\complexs$ and
its complex conjugate given by complex conjugation.

Note that  $[\tilde R,\tilde U]=0$.
The bundle $\tilde \bV$ gives rise to a geometric family $\tilde \cH=\tilde\cH^+\cup_B(\tilde \cH^-)^{op}$, where the underlying fiber bundle of
$\tilde \cH^\pm$ is again $\id\colon B\rightarrow B$, and
the twisting bundles are $\tilde \bV^\pm$.
For each
$\alpha\in [0,\pi/2]$ the operator
$$\cos(\alpha) \tilde U+\sin(\alpha)\tilde R$$ defines a taming $\tilde \cH_{t_\alpha}$. The family $(\tilde \cH_{t_\alpha})_\alpha$ 
defines a taming $\hat\cH_t$ of $\hat \cH:=\pr^*_B \tilde \cH$
over $[0,\pi/2]\times B$.
A computation shows that
$d\eta(\hat \cH_t)=\pr^*_B \ch(\nabla^{\tilde V})=0$.
We conclude the following equality of de Rham cohomology classes
$$[\eta^{4k+1}(\tilde \cH_{t_0})]=[\eta^{4k+1}(\tilde \cH_{t_{\pi/2}})]\ .$$
An inspection of the definitions shows that
\begin{eqnarray*} 
\eta^{4k+1}(\tilde \cH_{t_0})&=&2\eta^{4k+1}(\cH_t)\\
\eta^{4k+1}(\tilde \cH_{t_{\pi/2}})&=&(X^*\circ i^*)\eta^{4k+1}(\cG_t).
\end{eqnarray*}
We conclude with Lemma \ref{coonex} that
$$\hat c_{4k+2}(\cH)=X^*[i^* \frac{(2k)!}{2} \eta^{4k+1}(\cG_t) ]_{\R/\Z}=X^* \bar d^\R_{4k+1}\ .$$
\end{proof}

\subsection{The topological interpretation of $d^{-n}_{B,4k+1-n}$}

\subsubsection{}
Recall that the classifying space of $KO^{-n}$ is $\Omega^nBO$.
In view of the fibration
$$\Omega^nU/O\rightarrow \Omega^nBO\rightarrow \Omega^nBU$$
we see that the classifying map $X\colon B\rightarrow \Omega^n BO$ of an element
$x\in U^{-n}_\infty(B)$ factors (up to homotopy) over $\Omega^nU/O$,
since then the
composition $B\xrightarrow{X} \Omega^n BO\to\Omega^n BU$ is null homotopic.

%Applying the loop map of Definition
%\ref{def:loop_transgression} $n$ times to  $d_{4k+1}\in
%H^{4k+1}(U/O,\Q)$ we obtain a class
%$\Omega^n d_{4k+1}\in H^{4k+1-n}(\Omega^n U/O,\Q)$.
Let $\Omega^n:H^{4k+1}(U/O,\R)\rightarrow H^{4k+1-n}(\Omega^n U/O,\R)$
be the $n$-fold iteration of the loop map introduced in Definition \ref{def:loop_transgression}.
\begin{theorem}\label{theo:top2}
We have
$$d^{-n}_{B,4k+1-n}(x)=X^* [\frac{ (m-1)!}{(2k)!}\Omega^n
d^\R_{4k+1}]_{\R/\Z}= X^*[\frac{(m-1)!}{(2k)!}\Omega^n
T(c_{4k+2}^{\reals})]_{\reals/\integers}\ , $$
where $m$ is determined by $2m=4k+3-n$ or $2m=4k+2-n$.
\end{theorem}
\begin{proof}
We employ Corollary \ref{jh}. Let $\hat x\in \widetilde{KO}^0(\Sigma^n
B)$
correspond to $x\in KO^{-n}(B)$ under the identification
$ KO^{-n}(B)\cong \widetilde{KO}^0(\Sigma^n
B)$.
Let $\tilde x\in \widetilde{KO}^0(S^n\times B)$ be the pull-back of
$\hat x$ under the natural
map $S^n\times B\rightarrow \Sigma^n B$.

Note that the classifying map $Y\colon \Sigma^n B\to
U/O$ of $\hat x$ is the adjoint of $X:B\rightarrow \Omega^nU/O$, and that the composition $\tilde X\colon S^n\times B\rightarrow U/O$ of  the projection
$S^n\times B\to \Sigma^n B$ and $Y$ is  the classifying map
of $\tilde x$.

% and $Y\colon \Sigma^n B\to U/O$ the classifying map of $x$.
%Note that $\tilde X$ is the composition of $Y$ and the projection
%$S^n\times B\to \Sigma^n B$. Moreover, $Y\colon \Sigma^n B\to U/O$ and
%$X\colon B\to \Omega^n U/O$ are adjoint to each other.

Then we have
$$d^{-n}_{B,4k+1-n}(x)=[(m-1)!\int_{S^n\times B/B}\eta^{4k+1}(\cH_t)]_{\R/\Z}\ ,$$
where $\cH_t$ is constructed as in the proof of Theorem \ref{top1}.
In that proof we have also shown that
\begin{equation*}
[\eta^{4k+1}(\cH_t)]=[\frac12 (\tilde X^*\circ
i^*)\eta^{4k+1}(\cG_t)]\ .
\end{equation*}
We now apply Lemma \ref{coonex} in order to conclude that
\begin{equation*}
[\eta^{4k+1}(\cH_t)]=\frac{\tilde X^* d^\R_{4k+1}}{(2k)!}\ . 
\end{equation*}
Thus
$$d^{-n}_{B,4k+1-n}(x)=[\frac{(m-1)!}{(2k)!} \int_{S^n\times B/B } \tilde X^* d^\R_{4k+1}]_{\R/\Z}\ .$$
The assertion now follows from the general fact that for any $z\in H^{4k+1}(U/O,\R)$ we have 
$$\int_{S^n\times B/B } \tilde X^* z= \Sigma^n Y^*z =X^*\Omega^n z\ ,$$
where $\Sigma$ is the suspension isomorphism. For the first equality we use that
integration over the fiber essentially is the suspension isomorphism
in the above construction. The second equality is a special case of the
relation between suspension and
loop homomorphism proved in  Lemma \ref{lem:looping_suspension_and pull
  back}. 
\end{proof} 

\subsubsection{}
 \label{rem:special implies general}
 Theorem \ref{top1} and Theorem \ref{theo:top2} 
give a topological description of the value of 
 $d^{-n}_{B,4k+1-n}(x)$ only under the additional assumption that  $x\in
 U^{-n}_\infty(B)\subset 
 U^{-n}_{4k+2-n}(B)$. In order to see that this 
 determines $d^{-n}_{B,4k+1-n}$ completely we argue as follows.

Let $x\in U^{-n}_{4k+2-n}(B)$. Then the cohomology class
  $d^{-n}_{B,4k+1-n}(x)$ of degree $4k+1-n$ is
  determined by its restriction $d^{-n}_{B,4k+1-n}(x)_{|B^{4k+1-n}}$
    to the $(4k+1-n)$-skeleton
  $B^{4k+1-n}$ of $B$. We have
$x_{|B^{4k+1-n}}\in U^{-n}_\infty (B^{4k+1-n})$. Thus we
know the topological description of 
 $d^{-n}_{B^{4k+1-n},4k+1-n}(x_{|B^{4k+1-n}})$, which is equal to
$d^{-n}_{B,4k+1-n}(x)_{|B^{4k+1-n}}$ by naturality. 

\subsection{Explicit calculation of the universal class}
\label{sec:transgression}

Theorem \ref{theo:top2} does give a topological interpretation of our
invariant \ref{def:of natural transformation}. However,  we want to be
even more precise and explicitly compute the corresponding universal
cohomology class
$$[\frac{(m-1)!}{(2k)!}\Omega^n
d_{4k+1}]_{\Q/\Z}=[\frac{(m-1)!}{(2k)!}\Omega^n
T(c^\Q_{4k+2})]_{\Q/\Z}\in H^{4k+1-n}(\Omega^n U/O,\Q/\Z)\ ,$$
where $m$ is determined by $2m=4k+3-n$ or $2m=4k+2-n$.
In particular, we will show that for
half of the parameters $n$ (mod $(8)$) this class vanishes

We will make use of many of the results about the cohomology of $BO$,
$BU$ and their loop spaces collected in Appendix
\ref{sec:cohomology-bo-bu}.

\subsubsection{}
Consider the map of fibrations
\begin{equation*}
  \begin{CD}
    U @>>> EU @>>> BU\\
    @VV{p}V @VVV @VV{\id}V\\
    U/O @>i>> BO @>{Bc}>> BU,
  \end{CD}
\end{equation*}
where the upper row  is the universal principal $U$-bundle, and the
lower row  is obtained from the upper by dividing out the subgroup
$O$. 

\subsubsection{}
By Lemma
\ref{lem:transgressive classes}, 
$c_{4k+2}^\rationals \in H^{4k+2}(BU,\rationals)$ is transgressive in the second (and
of course also in the first) fibration. We have to compute the
transgression $T(c_{4k+2}) \in H^{4k+1}(U/O,\rationals)$. To do this,
we observe that the upper fibration is 
%by Bott periodicity 
the path
space fibration, and therefore by Lemma \ref{lem:loop is
  transgression} the transgression $T_U$ of this
fibration coincides with the
loop homomorphism. By Theorem \ref{theo:loop map for BO etc}  we
obtain $T_U(c_{4k+2}) = c_{4k+1}$ (even in
integral cohomology). Moreover, transgression is natural, therefore
\begin{equation*}
  p^*(T(c^\rationals_{4k+2})) = T_U(c^\rationals_{4k+2}) =
  c^\rationals_{4k+1} \in H^*(U,\rationals).
\end{equation*}
By Theorem \ref{theo:cohom of BO etc}, Theorem \ref{theo:complex Bott}
and Table \ref{table:maps-between-loop their effect in homology},
$p^*\colon H^*(U/O,\rationals)\to H^*(U,\rationals)$ is injective and
$p^*(a_{4k+1}/2) = c_{4k+1}$. 

\begin{notation}\label{not6}
In order to avoid an inflationary appearance of the exponent ${}^\Q$ 
  from now on we will use the same symbol for an integral cohomology
  class and its image in rational cohomology. It will be clear from
  the context which meaning the symbol has.
\end{notation}

Consequently (with the new convention \ref{not6} ) we can write
\begin{equation*}
  T(c_{4k+2}) = \frac{1}{2} a_{4k+1} \in H^{4k+1}(U/O;\rationals).
\end{equation*}

\subsubsection{}

Our next goal is the calculation of  $$\Omega^n T(
c_{4k+2})=\Omega^n( \frac{1}{2} a_{4k+1}) \in
H^{4k+1-n}(\Omega^nU/O,\rationals)\ .$$ 
We consider the fibration
 \begin{equation}\label{nlo}
  \Omega^nU/O\to \Omega^n BO\to \Omega^n BU
\end{equation}
which is the $n$-fold loop of the fibration considered above. 
%By Proposition \ref{prop:transgression and loop} 
%transgression commutes with the loop map.
%Therefore 
%$$T\Omega^n (c_{4k+2})=\Omega^n( \frac{1}{2} a_{4k+1})\ ,$$
%where the transgression is associated to the fibration
%(\ref{nlo}), and it suffices to compute $T\Omega^n(c_{4k+2})$.

%$\Omega^n(T
%c_{4k+2})= T\Omega^n(c_{4k+2})$).
\subsubsection{}
In the following, we use the Bott periodicity maps to identify
$\Omega^nU/O$ with the spaces listed in Theorem \ref{theo:Bott periodicity list}:

\begin{tabular}[l]{l|lllllllllllll}
$n$ & 0 & 1 & 2 & 3 & 4 & 5 & 6 & 7\\ \hline
  $\Omega^nU/O$ & U/O &$ BO\times\integers$ & O & O/U & U/Sp &
  $BSp\times\integers$ & Sp & Sp/U 
\end{tabular}

\subsubsection{}
%Since transgression commutes with applying the functor $\Omega$, we
%simply have to compute $\Omega^n(\frac{1}{2} a_{4k+1}) \in
%H^{4k+1-n}(\Omega^n(U/O);\rationals)$.
Unfortunately, our knowledge
about the map $\Omega^n$ is not complete enough to calculate
$\Omega^n (Tc_{4k+2})$ directly. 
We use the following trick:

Using that map $\alpha:U/O\to U$ (compare Subsection
\ref{sec:maps-between-our}, Item \ref{item:U_modO_to_U}) we have 
$\frac{1}{2}a_{4k+1}=\alpha^*
\frac{1}{2}c_{4k+1}$.  
Therefore, 
%$T(\Omega^n
%c_{4k+2})=\Omega^n T(c_{4k+2})=
$$\Omega^n(\frac{1}{2} a_{4k+1})=(\Omega^n\alpha)^*
\Omega^n(\frac{1}{2}c_{4k+1})
\ .$$
%\in H^{4k+1-n}(\Omega^n
%U,\rationals)$.
\subsubsection{}
We shall first compute $\Omega^n(c_{4k+1})$. By Theorem \ref{theo:loop map for BO etc} 
\begin{equation*}
\Omega(c_{4k+1}) = (2k)! \ch_{4k}.
\end{equation*}
Note that $\Omega
BU=\Omega(BU\times\integers)$ so that we can iterate the argument.

Next, $(2k)!\ch_{4k}= 2k \cdot c_{4k} + \text{\em decomposable}$.  
We conclude that $$\Omega( (2k)! \ch_{4k})
= 2k \Omega(c_{4k}) = 2k c_{4k-1}\ .$$ 

\subsubsection{}
Now, an easy induction allows us
to compute $\Omega^n(c_{4k+1})$ for each $n\in\naturals$. However,
because of the factors appearing in our formulas, we really have to
study $\frac{(m-1)!}{(2k)!} \Omega^{n}(c_{4k+1})$, where $m=2k+1-\frac{n}{2}$
if $n$ is even, and $m=2k+1 -\frac{n-1}{2}$ if $n$ is odd. In the induction,
this factor cancels the factors (like $(2k)$) which show up in the
calculations above, and we get (with $m$ depending on $n$ and $k$ as above)
\begin{equation}\label{eq:Omega c}
  \Omega^n( \frac{(m-1)!}{(2k)!} c_{4k+1}) =
  \begin{cases}
    c_{4k+1 -n} & n\text{ even}\\
    (2k -\frac{n-1}{2})!\ch_{4k +1 - n} & n\text{ odd} \ .
  \end{cases}
\end{equation}

\subsubsection{}
From this and the calculation of the map in cohomology induced by 
$\Omega^n\alpha:\Omega^nU/O\to
\Omega^nU$, we read off the cohomology classes we are
interested in. Since we are really only interested in the image of the
 class under the map in cohomology induced by the map of
coefficients $\rationals\to\rationals/\integers$, we obtain the following
list:
\begin{theorem}\label{theo:final answaer} 
\begin{equation}\label{eq:final result for d}
  \begin{split}
    [\frac{(m-1)!}{(2k)!} T(c_{4k+2})]_{\Q/\Z} &= [\frac{1}{2} a_{4k+1}]_{\Q/\Z} \in
    H^{4k+1}(U/O;\rationals/\integers)\\
    [\Omega(\frac{(m-1)!}{(2k)!} T(c_{4k+2})) ]_{\Q/\Z} &=
    \begin{cases}
      [\frac{1}{2} \ch_{0}]_{\Q/\Z} \in H^0(BO\times\integers,\rationals/\integers)  & 4k+1=n\\
      0 \in H^{4k}(BO\times\integers; \rationals/\integers) & 4k+1>n
      \end{cases}\\
      [\Omega^2(\frac{(m-1)!}{(2k)!}  T(c_{4k+2})) ]_{\Q/\Z} &= 0 \in H^{4k-1}(O,\rationals/\integers)\\
      [\Omega^3(\frac{(m-1)!}{(2k)!}  T(c_{4k+2})) ]_{\Q/\Z} &= 0 \in
      H^{4k-2}(O/U,\rationals/\integers)\\
      [\Omega^4(\frac{(m-1)!}{(2k)!}  T(c_{4k+2})) ]_{\Q/\Z} &= 0 \in
      H^{4k-3}(U/Sp,\rationals/\integers)\\
      [\Omega^5(\frac{(m-1)!}{(2k)!}  T(c_{4k+2})) ]_{\Q/\Z}  &= 0 \in
      H^{4k-4}(BSp\times\integers,\rationals/\integers)\\
      [\Omega^6(\frac{(m-1)!}{(2k)!}  T(c_{4k+2})) ]_{\Q/\Z} &=
      [\frac{1}{2}y_{4k-5}]_{\Q/\Z} \in
      H^{4k-5}(Sp,\rationals/\integers)\\
      [\Omega^7(\frac{(m-1)!}{(2k)!}  T(c_{4k+2})) ]_{\Q/\Z} &=
      [\frac{1}{2}(c_{4k-6}+ c_2c_{4k-8}+\dots+c_{2k-4} c_{2k-2}) ]_{\Q/\Z}\\
            & \qquad\in
      H^{4k-6}(Sp/U,\rationals/\integers)\\
    \end{split}    
\end{equation}
For $n\ge 8$, the answer can be read off from the list by reduction
mod $8$ by Bott periodicity.

In particular, 
the natural transformation $d^{-n}_{B,4k+1-n}$ vanishes for $n$
congruent to
$2,3,4,5$ mod $8$, and,  if $k>0$, also for $n\equiv 1 \pmod 8$.
In the other cases, since the universal classes are non-trivial, there
are non-trivial examples.
\end{theorem}
%\subsubsection{}
%Note that the statement that the {\em natural transformation}
%$d^{-n}_{B,4k+1-n}$ is non-trivial, means that there exist spaces
%$B$ for which this map is non-trivial.
\subsubsection{}
\begin{proof}
  As observed above, we simply have to take the cohomology
  classes on the right-hand side of  equation \eqref{eq:Omega c},
  divide them by $2$, and then
  apply the map from rational cohomology to cohomology with coefficients
  in $\rationals/\integers$. Finally through $\Omega^n\alpha$  we pull
  back the result to $\Omega^n U/O$. In this step we use the
  results of Subsection \ref{sec:maps-between-our}. 

\subsubsection{}
  Note first that $x_n:=c_{4k+1-n}$ ($n$ even), and respectively, $x_n:=(2k
  -\frac{n-1}{2})!\ch_{4k +1 - n}$ ($n$ odd) belong to the integral
  lattice in rational cohomology. 
  Therefore $[\frac{1}{2}(\Omega^n\alpha)^*\Omega^n x_n]_{\R/\Z}=0$
if $\frac{1}{2}(\Omega^n\alpha)^*$ preserves the integral lattices.
This is the case 
%The pull back of $1/2$ these classes from $H^*(\Omega^n
%  U;\rationals/\integers)$ 
%  to $H^*(\Omega^n U/O;\rationals/\integers)$ therefore is annihilated
%  by $2$, and vanishes
 whenever  $\Omega^n\alpha$
   maps the Chern classes to twice a generator of the
  integral cohomology, i.e. if $\Omega^nU/O$ equals $O$, $O/U$, or
  $U/Sp$
by the table \ref{table:maps-between-loop their effect in homology}.
This observation accounts for the zeros for $n=2,3,4$ in the theorem.

\subsubsection{}
Because of Proposition \ref{prop:transgression and loop}
$$
  [T(\Omega^5(\frac{(m-1)!}{(2k)!}  c_{4k+2}))]_{\Q/\Z}=
%\in 
%      H^{4k-4}(BSp\times\integers;\rationals/\integers)$ equals
      [\Omega T(\Omega^4(\frac{(m-1)!}{(2k)!}  c_{4k+2}))]_{\Q/\Z} =0\
      ,$$ we obtain the zero
      for $n=5$.
\subsubsection{}
We now discuss the case $n=1$.
We have $[T(c_{4k+2})]_{\Q/\Z}=l_*(w_{2k}\cup
      w_{2k+1})$, where $l_*$ is induced
      by the map of coefficients
      $\integers/2\integers \to \rationals/2\integers
      \xrightarrow{\cdot 1/2} \rationals/\integers$
(compare \ref{lsdeef} and Corollary \ref{spi1}).
%In fact (here $a_{4k+1}$ is the integral class, compare Theorem
%\ref{theo:cohom of BO etc}, c.)
%since $[a_{4k+1}]_{\Z/2\Z}=w_{2k}\cup w_{2k+1}+\beta( \ldots)
%      \in H^{4k+1}(U/O,\integers/2\integers)$, and the image of
%      the Bockstein operator $\beta(\ldots)$ gets annihilated  when mapped further to
%      $H^{4k+1}(U/O),\rationals/2\integers)$. 
If $k>0$, then
      $\Omega(w_{2k}\cup w_{2k+1})=0$ since the loop
      map is applied to a decomposable class. 
Note that $\frac{(m-1)!}{(2k)!}=1$ in this case.
Thus $$   [\Omega (\frac{(m-1)!}{(2k)!}T(
      c_{4k+2}))]_{\Q/\Z}=l_*\Omega(w_{2k}\cup w_{2k+1}) =0$$ for
      $k>0$.

%\in
%      H^{4k}(BO\times\integers,\rationals/\integers)$ is obtained from
%      $\Omega(w_{2k}\cup w_{2k+1}) \in
%      H^{4k}(BO\times\integers,\integers/2\integers)|$ by mapping to
%      $H^{4k}(BO\times\integers,\rationals/\integers)$ as above, and
%      therefore vanishes for $k>0$.
\subsubsection{}
      For the calculation of $[\Omega^7(\frac{(m-1)!}{(2k)!}
      T(c_{4k+2}))]_{\Q/\Z}$ we proceed as follows.
      The class $(2k-3)!\ch_{4k-6}\in H^{4k-6}(Sp/U,\Q)$ belongs to
      the integral lattice. In fact, if we write 
$(2k-3)!\ch_{4k-6}=\sigma_{4k-6}(c_2,c_4,\dots,c_{4k-6})$ with the
Newton polynomial $\sigma_{4k-6}$,
then the right-hand side can be interpreted as an integral cohomology
class in $H^{4k-6}(Sp/U,\Z)$.
We now have $$l_*[\sigma_{4k-6}(c_2,c_4,\dots,c_{4k-6})]_{\Z/2\Z}=[\Omega^7(\frac{(m-1)!}{(2k)!}
      T(c_{4k+2}))]_{\Q/\Z}\ .$$  
The cohomology $H^{4k-6}(Sp/U,\Z/2\Z)$ is an exterior algebra
generated by  $[c_2]_{\Z/2\Z},[c_4]_{Z/2\Z},\dots$.
Considered ---in this algebra--- the Newton polynomial 
\begin{eqnarray*}
\lefteqn{\sigma_{4k-6}(c_2,c_4,\dots,c_{4k-6}) }&&\\&=&-(2k-3)
        \sum_{i_1+2i_2+\dots (2k-3)i_{2k-3}={2k-3}} (-1)^{i_1+\dots+i_{2k-3}}
        \frac{(i_1+\dots+i_{2k-3}-1)!}{i_1!\dots i_{2k-3}!}
        c_2^{i_1}\dots c_{4k-6}^{i_{2k-3}}.\end{eqnarray*}
It simplifies considerably and
%we use again the fact that we can first look at the
%      integral class $2T(\Omega^7(\frac{(m-1)!}{(2k)!}  c_{4k+2}))\in
%      H^*(Sp/U;\integers)$, map to $H^*(Sp/U;\integers/2\integers)$
%      and then to $H^*(Sp/U;\rationals/\integers)$ as above. Since the
%      Chern character is given by the Newton polynomial and the Chern
%      classes are pulled back to the classes $c_*$ by
%      \ref{table:maps-between-loop their effect in homology}, this is
%      simply given by the Newton polynomial
    %  \begin{multline*} 2T(\Omega^7(\frac{(m-1)!}{(2k)!}  c_{4k+2})) =
%        \sigma_{4k-6}(c_2,c_4,\dots,c_{4k-6}) =\\ -(2k-3)
%        \sum_{i_1+2i_2+\dots (2k-3)i_{2k-3}={2k-3}} (-1)^{i_1+\dots+i_{2k-3}}
%        \frac{(i_1+\dots+i_{2k-3}-1)!}{i_1!\dots i_{2k-3}!}
%        c_2^{i_1}\dots c_{4k-6}^{i_{2k-3}}.
%        \end{multline*}
%        In the exterior algebra $H^*(Sp/U;\integers/2\integers)$
%        generated by $c_2,c_4,\dots$ this simplifies considerably and
%        gives (modulo $2\integers$) 
gives exactly the expression asserted, if we
use that 
$$l_*[c_{2l}]_{\Z/2\Z}=[\frac{1}{2} c_{2l}]_{\Q/\Z}\ .$$
\end{proof}  

\subsection{The relation with ordinary characteristic classes}\label{cl7}
\subsubsection{}
%\begin{remark}
Let us consider the fibration
$$\Omega^n U/O\stackrel{\Omega^ni}{\rightarrow} \Omega^n
BO\rightarrow \Omega^n BU\ .$$ 
We have constructed and calculated the universal cohomology class $$[\Omega^n(\frac{(m-1)!}{(2k)!}  T(c_{4k+2})) ]_{\Q/\Z}\in
H^{4k+1-n}(\Omega^n U/O,\Q/\Z)\ .$$
If this class would be of the form
$(\Omega^ni)^* u$ for some $u\in H^{4k+1-n}(\Omega^nBO,\Q/\Z)$,
then the invariant $d^{-n}_{B,4k+1-n}(x)$, $x\in U^{-n}_{4k+1-n}(B)$,
could be expressed in terms of familiar characteristic classes of the
element $x\in KO^{-n}(B)$. 

\subsubsection{}
In the case $n=0$ we indeed have
$$[\Omega^n(\frac{(m-1)!}{(2k)!}  T(c_{4k+2}))
]_{\Q/\Z}=[T(c_{4k+2})]_{\Q/\Z}=i^*l_*(w_{2k}\cup w_{2k+1})\ $$
(above we have written $w_l$ for $i^*w_l$ in order to save
notation, but in the current discussion it  makes sense not to omit $i^*$).
In particular we can extend $d^0_{B,4k+1}$ to all of $KO^0(B)$
by setting $d^0_{B,4k+1}(x):=l_*(w_{2k}\cup w_{2k+1})$.

\subsubsection{}
In the case $n=1$ and $4k=0$ it is obvious that the class
comes from $\Omega^1 BO$. As we have seen in example \ref{lmk},
if $x \in KO^{-1}(X)$ is represented by a family of anti-selfadjoint
real Fredholm operators $(F_b)_{b\in B}$, then $d^0_{B,0}(x)$ is
represented by
the locally constant $\R/\Z$-valued function $b\mapsto [\frac{1}{2}\dim\ker F_b]_{\R/\Z}$

\subsubsection{}

%This means that our classes are not even special
%  to the situation that the complexification is trivial, but are
%  defined for arbitrary real families.

%  The same is true for $n=1$. This is a rather trivial result, since
%  the only non-vanishing class exists in
%  $H^0(BO\times\integers;\rationals/\integers)$.

In the case $n=6$ the class $$[\Omega^6(\frac{(m-1)!}{(2k)!}  T(c_{4k+2}))
]_{\Q/\Z}=[\Omega^6(\frac{1}{2k(2k-1)(2k-2)}  T(c_{4k+2}))
]_{\Q/\Z}$$ definitely is not a pull-back from $\Omega^6BO=Sp/U$.
In fact
    $H^*(Sp/U,\rationals/\integers)$ is
  concentrated in even degrees, while our class is of odd degree.
We see that in this case our invariant $d^{-6}_{B,4k-5}$ is more exotic
and therefore more interesting. Unfortunately, we haven't been able to
produce simple examples of non-triviality for this invariant in Section
\ref{trio}. 
   
\subsubsection{}\label{sec:omega7i}

  For $n=7$, the map $Sp/U\to U/O$ induces a surjection in
  cohomology with $\integers/2\integers$-coefficients by Lemma
  \ref{lem:surjecti}. Therefore, all our
  classes pull back from $U/O$.

\subsection{Extendibility}

\subsubsection{}
Given $x\in KO^{-n}(B)$, in order to define $d^{-n}_{4k+1-n}(x)$ using
topology we had to assume that  $x\in U^{-n}_\infty(B)$.
Our analytic definition however works under the weaker condition that
$x\in U_{4k+2-n}^{-n}(B)$. Of course, if $B^{4k+1-n}\subset B$ denotes
a $4k+1-n$-skeleton, we have
$x_{|B^{4k+1-n}}\in U^{-n}_\infty(B^{4k+1-n})$. We also have seen that
$d^{-n}_{4k+1-n}(x_{|B^{4k+1-n}})=d^{-n}_{4k+1-n}(x)_{|B^{4k+1-n}}$
determines $d^{-n}_{4k+1-n}(x)$ uniquely.
The interesting feature of the analytic definition is
that it shows that $d^{-n}_{4k+1-n}(x_{|B^{4k+1-n}})$
admits an extension from $B^{4k+1-n}$ to $B$.

In the following Lemmas we give an alternative proof of this property.
\subsubsection{}

Let $f\colon B\to \Omega^n BO$ be a map. Assume that the restriction
  $f^{r}:=f_{|B^r}\colon B^{r}\to \Omega^n BO$ of $f$ to a $r$-skeleton
  $B^r$ of
  $B$ factors over a map $g^{r}\colon B^{r}\to \Omega^nU/O$
  (i.e.~$f^{r}= \Omega^n i\circ g^{r}$, where $i\colon U/O\to BO$
  is as above). 
Assume further  that $r+n\equiv 1 \pmod{8}$. Let $R$ be some  abelian group.
\begin{lemma}\label{lem:extension of classes}
  If $x\in H^k(\Omega^n U/O,R)$ satisfies $2x=0$,  then the class 
$(g^{r})^*(x)$ extends from $B^r$
  to $B$.
\end{lemma}
 Since the map $H^r(B;R)\to H^r(B^{r};R)$ is 
  injective, this extension is unique. Note that the Lemma in
  particular applies
to the cohomology classes listed in Theorem
  \ref{theo:final answaer}. 
\subsubsection{}
\begin{proof}
Let  $K(R;r)$ denote the Eilenberg-Mac Lane space which represents
the functor $H^r(\dots,R)$. We represent the cohomology class $x$ by a
map
 $x\colon \Omega^nU/O\to K(R;r)$.  It suffices to show that
 $(g^{r})^*(x)$ extends to 
a $r+1$-skeleton $B^{r+1}$ of $B$ (such that $B^r\subset B^{r+1}$). In
fact, it then further extends to $B$ since
the inclusion $B^{r+1}\to B$ is a $r+1$-equivalence.

\subsubsection{}
  The universal example is given by the space $B=K^{r+1}$ which is
  obtained from $\Omega^nU/O$ by attaching $r+1$-cells in such a
  way as to kill the kernel of $(\Omega^ni)_r:\pi_k\Omega^n(U/O)\to \pi_r(\Omega^n
  BO)$. Here $f$ is obtained from $\Omega^ni\colon
  \Omega^n(U/O)\to \Omega^n BO$, which extends to some map $f\colon
  K^{r+1}\to \Omega^n BO$ by the
  construction of $K^{r+1}$ (and elementary obstruction theory), and
  $g^r$ is the inclusion of the $r$-skeleton of $\Omega^n(U/O)$ (and
  therefore of $K^{r+1}$) into $\Omega^n(U/O)$.

In
  our case we have $\pi_r(\Omega^nU/O)\iso\integers$ and
  $\ker(\Omega^ni)_r\cong 2\Z$ (as follows from Bott periodicity and the long
  exact homotopy sequence (in low degrees) of $U/O\to BO\to BU$). If
  $\phi\colon S^r\to \Omega^n(U/O)$ represents a generator of
  $\pi_r(\Omega^nU/O)$, and if $h\colon
  S^r\to \Omega^nU/O$ represents twice this generator, i.e.~a generator of
  $\ker(\Omega^ni)_r$,   then $h^*
  (g^{r})^* x = 2 \phi^* (g^{r})^* x=0$. 
Thus the map $x\circ h$ is null
  homotopic, and therefore $(g^{r})^*x$ extends to $K^{r+1}$. \end{proof}

\subsubsection{}
Let still $f\colon B\to \Omega^n BO$ be a map and assume that the restriction
  $f^{4k+1-n}:=f_{|B^{4k+1-n}}\colon B^{4k+1-n}\to \Omega^n BO$ of $f$
  to a $4k+1-n$-skeleton
  $B^{4k+1-n}$ of
  $B$ factors over a map $g^{4k+1-n}\colon B^{4k+1-n}\to \Omega^nU/O$
  (i.e.~$f^{4k+1-n}= \Omega^n i\circ g^{4k+1-n}$, where $i\colon U/O\to BO$
  is as above).

\begin{lemma}\label{lem:extend}
  If $x=[\Omega^n(\frac{(m-1)!}{(2k)!}T(c_{4k+2}))]_{\Q/\Z}\in
  H^{4k-n+1}(\Omega^nU/O)$ is one of the classes of Theorem
  \ref{theo:final answaer}, then $(f^{4k+1-n})^*x$ extends from
  $B^{4k+1-n}$ to $B$.
\end{lemma}

\subsubsection{}\label{sec:half_classes}
Note that half of the cases are already covered by Lemma
\ref{lem:extension of classes}, namely whenever the dimension
condition is satisfied, i.e.~when $(4k+1-n)+n \equiv 1\pmod 8$, in
other words, if $k$ is even.

\subsubsection{}

Moreover, the cases $n\equiv 0\pmod 8$ as well as $n\equiv 1 \pmod 8$ (and
$4k+1=n$) are
trivial, because in these cases we have seen that the characteristic
class $x$ already pulls back from $\Omega^n BO$ to $\Omega^n U/O$: it
is expressed in terms of Stiefel-Whitney classes in the first case,
and in terms of the dimension of the bundle in the second case.

\subsubsection{}

We use the proof Lemma \ref{lem:extension of classes} to deal with the
remaining cases. This proof shows that it suffices to treat the case
$B^{4k+1-n}=S^{4k+1-n}$ and to show that the pullback class
$(f^{4k+1-n})^*x$ 
vanishes for arbitrary $f\colon S^{4k-n+1}\to \Omega^nU/O$ (and therefore extends over the disc $D^{4k+2-n}$).

\subsubsection{}
\label{sec:chern_integral}

Observe that, by Equation
  \eqref{eq:Omega c}, the cohomology class $x$ is obtained as pull
  back of $\frac{1}{2}c_{4k+1-n}$ or
  $\frac{1}{2}(2k-(n-1)/2)! \ch_{[4k-(n-1)]}$ from $U$ or $BU$,
  respectively (depending on the parity of $n$).
However, on all spheres the Chern character is
  integral, i.e.~for an arbitrary map $f\colon S^k\to BU$,
  $f^* \ch \in H^*(S^k;\integers)$.

\subsubsection{}
  
  If
  $n\equiv 7 \pmod 8$ and $k$ odd (and $4k>n$) then $4k+1-n\ge 4$. This implies
  $\frac{1}{2}(2k-(n-1)/2)!\in\integers$. Therefore by \ref{sec:chern_integral}
  the cohomology class  $\frac{1}{2}(2k-(n-1)/2)! \ch_{[4k-(n-1)]}$
  pulls back to $0$ in $H^*(S^{4k-(n-1)};\rationals/\integers)$ for an
  arbitrary map $S^{4k-(n-1)}\to BU$. As observed in
  \ref{sec:half_classes}, $n\equiv 7 \pmod 8$ and
  $k$ even is covered by Lemma \ref{lem:extension of classes}.

\subsubsection{}

  For $n\equiv 6\pmod 8$ and an arbitrary map $f\colon S^{4k-n+1}\to U$, by Lemma
  \ref{lem:looping_suspension_and pull back} and Theorem
  \ref{theo:loop map for BO etc} 
  \begin{equation*}
    f^* (\frac{1}{2} c_{4k+1-n}) = \Sigma^{-1} F^*(\frac{1}{2} \Omega
    c_{4k+1-n}) = \Sigma^{-1}(F^*\frac{1}{2} (2k-\frac{n}{2})!
    \ch_{[4k-n]}).
  \end{equation*}
  Here, $\Sigma\colon H^{4k-n}(S^{4k-n})\to H^{4k-n+1}(S^{4k-n+1})$ is the suspension
  isomorphism and 
  \begin{equation*}
F\colon S^{4k-n}\to \Omega U=BU\times \integers
\end{equation*}
is
  the adjoint of $f\colon \Sigma S^{4k-n}=S^{4k-n+1}\to U$. Again, if
  $4k-n\ge 4$ then $(2k-n/2)!\in 2\Z$ and therefore by
  \ref{sec:chern_integral} $F^*(\frac{1}{2}(2k-n/2)!
  \ch_{[4k-(n-1)]}=0$. However, if $4k-n=2$ then, since $n\equiv
  6\pmod 8$ $k$ is even and therefore $(4k+1-n)+n\equiv 1\pmod 8$,
  such that this case is covered by Lemma \ref{lem:extension of
    classes}.

\subsubsection{}
If $n\equiv 2,3,4,5\pmod 8$, then $x=0$, which trivially extends. This
concludes the proof of Lemma \ref{lem:extend}. \hB

\begin{appendix}

\section{Transgression}
\label{sec:appendix_transgresseion}

\subsection{Transgression in cohomology}

\subsubsection{}

In this section, we want to recall the general definition of
transgression and  its basic properties. Special cases are
``suspension'' or ``looping''. All of this is well known, and included
here for the convenience of the reader.

\subsubsection{}

The situation is the following: let $f\colon E\to B$ be a map, and
$b\in B$ a point. Write $F_b:= f^{-1}(b)$. Let $i\colon F_b\to E$ be
the inclusion. Let $H^*$ be any (generalized) cohomology theory. In
the following, the loop spaces $\Omega B$ are defined with respect to
the basepoint $b$. Suspensions are reduced suspensions.

\subsubsection{}
The adjoint of the identity map $\Omega B\to\Omega B$ gives a
  canonical map $\Sigma\Omega B\to B$. This induces $H^*(B)\to
  H^*(\Sigma\Omega B)$. 
\begin{definition}\label{def:loop_transgression}
We define the loop map 
\begin{equation*}
    \Omega\colon H^*(B)\to H^{*-1}(\Omega B).
  \end{equation*} as the 
composition of $ H^*(B)\to
  H^*(\Sigma\Omega B)$ with the suspension isomorphism
  $H^*(\Sigma\Omega B)\rightarrow H^{*-1}(\Omega B)$ 
  \end{definition}
  By construction and functoriality of the suspension isomorphism,
  the loop map is functorial, too

\subsubsection{}\label{cofs}
Given the map $f\colon E\to B$, consider the cofibration sequence
$E\to Zf\to Cf$, where $Zf$ is the mapping cylinder and $Cf$ the
mapping cone. The inclusion $B\to Zf$ is a homotopy equivalence. The
long exact sequence in cohomology of this cofibration gives
\begin{equation*}
\begin{array}{ccccccccc}
&&&& H^k(B)&&&&\\
&&&&\|&&&&\\
  \cdots &\to& H^k(Cf)&\to& H^k(Zf)& \xrightarrow{f^*}&
  H^k(E)&\to&\cdots
\end{array}
\end{equation*}
In particular, $H^k(Cf)$ maps surjectively onto $\ker(f^*)\subset
H^k(B)$. 

\subsubsection{}\label{cofs2}
Consider now $\Omega Cf$. Since the composition
$F\xrightarrow{i} E\xrightarrow{f} B$ is the constant map to $b$, we can
define a canonical map
\begin{equation*}
l:   F\to \Omega Cf\ ,\quad  x\mapsto c_x,
\end{equation*}
where $c_x$ is the loop in $Cf$ with $c_x(0)=b\in B$, $c_x(t)=
(x,t)\in E\times (0,1)\subset Cf$, and $c_x(1)$ is the cone point in $Cf$.

Mapping $B$ to the second cone point gives the second map $j$ in the
cofibration sequence $B\to Cf\stackrel{j}{\rightarrow}\Sigma E$. From this we conclude that
the kernel of $H^k(Cf)\to H^k(B)$ equals $\im(j^*)$. The composition $F\to \Omega Cf\to \Omega \Sigma E$ can be
factored as $F\xrightarrow{i} E \to \Omega\Sigma E$, where the second
map is the adjoint of the identity map $\Sigma E\to\Sigma E$. 

\subsubsection{}
\begin{definition}
We define the transgression $T\colon H^k(B)\supset \ker(f^k)\to
  H^{k-1}(F)/\im(i^*)$ as the composition
  $$
    \ker(f^*)\iso H^k(Cf)/\im(j^*) \xrightarrow{\Omega}
    H^{k-1}(\Omega Cf)/ \im((\Omega j)^*)
    \to
    H^{k-1}(F)/\im(i^*)\ .
  $$\end{definition}
  Note that for the second map we used the factorization which shows
that $ \im((\Omega j)^*)$ goes into $\im (i^*)$.

%Alternatively, one can use the embedding $\Sigma B\into Cf$ and define
%the transgression as the composition
%\begin{multline*}
%    \ker(f^*)\iso H^k(Cf)/\im(H^k(\Sigma E)) \to H^k(\Sigma
%    F)/\im(H^k(\Sigma E)) \iso H^{k-1}(F)/\im(H^{k-1}(E)).
%\end{multline*}
%(We used the adjoint of $\Sigma B\into Cf$ to get $B\to \Omega Cf$).

It is clear from the construction that transgression is natural
with respect to the map $f:E\rightarrow B$, i.e.
given a diagram
$$\begin{array}{ccc}
E^\prime&\stackrel{H}{\rightarrow}&E\\
f^\prime\downarrow&&f\downarrow\\
B^\prime&\stackrel{h}{\rightarrow}&B
\end{array}$$
we have an equality of the form
$T^\prime\circ h^*=H^*\circ T$.

\subsubsection{}

\begin{definition}
The elements of $\ker(f^*)\subset H^*(B)$ are called 
 \emph{transgressive}.
 These are the classes whose
transgression is defined.
\end{definition}

\subsubsection{}
\begin{lemma}\label{lem:loop is transgression}
  The transgression in the fibration $\Omega B\to PB \to B$, where
  $PB$ is the (contractible) space of paths ending at $b$ coincides
  with the loop map.
\end{lemma}
\begin{proof}
  Carry out the construction. If $f\colon PB\to B$ is the start point
  projection, use the homotopy equivalence $Cf\to B$ which maps $ (p,s)\in PB\times
  (0,1)\subset Cf$ to $p(s)$ (recall that $s=1$ corresponds to the
  cone point). \end{proof}  
\subsubsection{}

\begin{lemma}\label{lem:sups_is_transgression}
  The transgression of $B\into CB\to \Sigma B$ is the suspension
  isomorphism $\Sigma\colon H^{k+1}(\Sigma B)\to H^k(B)$.
\end{lemma}
\begin{proof}
  Carry out the construction. Use the ``folding'' homotopy equivalence
  $Cf\to \Sigma B$ (where $f\colon CB\to \Sigma B$ is the
  projection). The composition of $\Sigma B\to  Cf$
  with this homotopy equivalence is the identity map. Starting with
  $H^*(\Sigma B)$, we have to pull back with this map and then use the
  suspension isomorphism (by naturality of the latter).
\end{proof}

 \subsubsection{}
  Let $f\colon \Sigma X\to Y$ be a map with adjoint $F\colon X\to
  \Omega Y$. Then we have a commutative diagram of
  fibrations
  \begin{equation*}
    \begin{CD}
      X\ @>>> CX@>>> \Sigma X\\
      @VV{F}V @VVV @VV{f}V\\
      \Omega Y @>>> PY @>>> Y  
    \end{CD}\quad .
  \end{equation*}
\begin{lemma}\label{lem:looping_suspension_and pull back}
  Then for each $a\in H^*(Y)$, 
  \begin{equation}
F^*(\Omega a) = \Sigma f^* a\ .
\end{equation}
%i.e.~the suspension isomorphism and loop homomorphism are intertwined
%by the map $f$ and its adjoint $F$.
\end{lemma}
\begin{proof}
  This follows from naturality of transgression and the fact that both
  $\Omega$ and $\Sigma$ are transgression homomorphisms by Lemma
  \ref{lem:loop is transgression} and Lemma \ref{lem:sups_is_transgression}.
\end{proof}  

\subsubsection{}
\label{ref:suspension_determines_transgression}
  In the construction of the transgression, we consider in particular
  the following commutative diagram of maps:
  \begin{equation}\label{eq:susp det trnas}
    \begin{CD}
      F @>{i}>> E @>{f}>> B\\
      @VV{\id}V @VVV @VVV\\
      F @>>> CE @>>> Cf\\
      @AA{\id}A @AAA @AAA\\
      F @>>> CF @>>> \Sigma F
    \end{CD}
  \end{equation}
  By naturality, the transgression homomorphism in $E \xrightarrow{f}
  B$ is determined by the transgression in $CF \to \Sigma F$ (this is
  of course, what we used in the construction), since $H^*(Cf)\to
  H^*(B)$ surjects onto the transgressive classes.

\subsubsection{}
\begin{lemma}\label{lem:transgression is compatible with natural transformations}
  Let $\Phi\colon H^*\to h^*$ be a natural transformation between
  generalized cohomology theories. Transgression commutes with this
  natural transformation.
\end{lemma}
\begin{proof}
  Let $f\colon E\to B$ be a continuous map. First observe that by
  naturality $\Phi$ maps $\ker(H^*(f))$ to $\ker(h^*(f))$ and
  $\im(H^*(i))$ to $\im(h^*(i))$, so that the assertion makes sense.

  The construction of the transgression homomorphism only uses maps
  induced from continuous maps between topological spaces (and their
  inverses) and the suspension isomorphism. By definition, a natural
  transformation between cohomology theories is compatible with such
  homomorphisms, and therefore also with transgression.
\end{proof}  
\subsubsection{}
 Let $F\xrightarrow{i} E\xrightarrow{f} B$ be a sequence of maps as
  above. This gives rise 
  to the transgression homomorphism $T\colon H^*(B)\supset \ker(f^*)
  \rightarrow 
  H^{*-1}(F)/\im(i^*)$.

  Applying the loop space functor we also get the sequence of maps
  \begin{equation*}
    \Omega F\xrightarrow{\Omega i} \Omega E\xrightarrow{\Omega f}
    \Omega B,
  \end{equation*}
  with associated transgression homomorphism
  \begin{equation*}
    T_\Omega\colon H^*(\Omega B)\supset \ker(\Omega f^*) \to
    H^{*-1}(\Omega F)/ \im((\Omega i)^*)\ .
  \end{equation*}
\begin{proposition}\label{prop:transgression and loop}
  If $x\in
  H^n(B)$ is transgressive, then
  \begin{equation*}
    \Omega T(x) = T_\Omega(\Omega x) \in H^{n-2}(\Omega F)/\im((\Omega i)^*).
  \end{equation*}
\end{proposition}
\begin{proof}
 We obtain the following commutative diagram
 \begin{equation}\label{eq:path space fibrations}
   \begin{CD}
    \Sigma \Omega F @>{\Sigma\Omega i}>> \Sigma\Omega E
    @>{\Sigma\Omega f}>> \Sigma\Omega B\\ 
     @VV{p_F}V @VVV @VV{p_B}V\\
     F @>{i}>> E @>{f}>> B
   \end{CD}
 \end{equation}
  where the vertical maps are adjoints of the identity maps $\Omega
  \cdot\to\Omega\cdot$. Since we work with the reduced suspension, the
  inclusion $\Sigma\Omega F\into \Sigma\Omega E$ is the fiber of
  $\Sigma\Omega f$. By naturality of the transgression, $p_F^*(Tx)=
  T_{\Sigma\Omega}(p_B^*x)$ for each transgressive class $x\in
  H^*(B)$. The suspension isomorphism maps by definition $p_F^*(Tx)$
  to $\Omega(Tx)$ and $p_B^*x$ to $\Omega x$. By Lemma
  \ref{lem:transgression is compatible with natural transformations} transgression
  commutes with the
  suspension isomorphism
(indeed the suspension isomorphism can be interpreted as a natural
  transformation between cohomology theories).  Therefore we have
  $
    \Omega(Tx) = T_\Omega(\Omega x)
$.\end{proof}

\subsection{Transgression and products}
\label{sec:transgr-prod}

\subsubsection{}
Let $\Delta\colon B\to B\times B$ be the
   diagonal map.  We still
   consider the map $f:E\rightarrow B$.
\begin{definition}
  \begin{enumerate}
\item A class $x\in H^*(B)$ is  called a \emph{non-trivial product}, if
  $x=\Delta^* y$ for some
   $y\in H^*(B\times B)$ such that 
   $(\id_B\times\{b\})^*y=0=(\{b\}\times\id_B)^*y$.  
\item
  We say that the first factor of a non-trivial product $x$ is
   \emph{transgressive}, if $x=\Delta^*y$ for an $y$ such that
   $(f\times \id_B)^* y = 0 \in H^*(E\times B)$, similarly we define
   that the second factor is transgressive.
\end{enumerate}
\end{definition}
\subsubsection{}
Note that if one of the factors of a product is transgressive,
then so is the product.
\begin{proposition}\label{prop:transgression_of_products}
  The transgression of a non-trivial product with at least one transgressive
  factor is zero.
\end{proposition}
\begin{proof}
  This follows from naturality of the transgression. Consider the diagram
  \begin{equation*}
    \begin{CD}
       F@>>> E @>>> B\\
      @VV{\id}V @VV{(\id_E\times f)\circ \Delta}V @VV{\Delta}V\\
      F @>>> E\times B @>>> B\times B\\
      @AA{\id}A @AA{\id_{B}\times \{b\}}A @AA{\id_{B}\times \{b\}}A\\
          F @>>> E @>>> B
 \end{CD}
\end{equation*}
Let us assume that $x$ is transgressive in the first factor.
By $T_i$ we denote the transgressions associated to the corresponding rows.
Then we have $T_1(x)=T_2(y)=T_3((\id_B\times\{b\})^*y)=0$, since 
$(\id_B\times \{b\})^* y =0$.
\end{proof}  

 \subsubsection{}
Let us consider the following example.
  Define $N:=T^2\setminus (D^2)^\circ$, i.e.~$N$ is the two torus with an open
  disc removed. Let $f\colon N\to T^2$ be the map which collapses the
  boundary of $N$ to one point. 
On the one hand,  the fundamental class $[T^2]\in
  H^2(T^2, \integers)$ is transgressive. On the other hand,  $[T^2]$
  is a non-trivial product of $1$-dimensional cohomology classes, and
  none of the factors is transgressive.

  Collapsing the complement of an open disc in $T^2$ to a point gives
  a degree $1$ map $g\colon T^2\to S^2$. If we write $S^2=\Sigma S^1$, we then
  get a diagram
  \begin{equation*}
    \begin{CD}
      S^1 @>>> CS^1=D^2 @>>> \Sigma S^1=S^2\\
      @AA{\id_{S^1}}A @AAA @AA{g}A\\
      S^1 @>>> N @>{f}>> T^2
    \end{CD}
  \end{equation*}
  Here, $g^*[S^2]=[T^2]$, where $[S^2]\in H^2(S^2,\integers)$ is the
  fundamental class. By naturality, $T([T^2])=T([S^2])= [S^1]\in
  H^1(S^1,\integers)$ is the fundamental class of $S^1$, in particular
  non-zero.

  This shows that at least one of the factors in Proposition
  \ref{prop:transgression_of_products} has to be transgressive for the
  assertion to hold.

\subsection[Transgression and relation with Bockstein]{Transgression in ordinary cohomology  and the relation
  with the Bockstein}
\label{sec:expl-cycl-transgr}
\subsubsection{}
We now want to describe how one can construct the
transgression in ordinary singular cohomology with coefficients on the
level of chains.
Let $f\colon E\to B$ be a map with fiber $i:F\hookrightarrow E$ over
$b\in B$. Let $R$ be an
  abelian group. Assume that $x\in H^k(B,R)$ is transgressive,
  i.e.  we have $f^*x=0$. We
choose a cocycle $c\in C^k(B,R)$ representing $x$. Then the cocycle
  $f^*c$ is a boundary, i.e. there exists a chain $c_0\in C^{k-1}(E,R)$ with
  $dc_0=f^*c$. The restriction of $c_0$ to $F$ is closed, since
$di^*c_0=i^*dc_0=i^*p^*c$.
  It follows that
  $i^*c_0$ represents a cohomology class $[i^*c_0]\in H^{k-1}(F,R)$. The
  cocycle $c_0$ is well defined only up to closed cocycles in $E$. It follows
  that the class $[c_0]$ is well defined only up to the image of
  $i^*$. Hence we get a well-defined class $\tilde T(x)\in
  H^{k-1}(F)/\im(i^*)$.
\subsubsection{}
\begin{proposition}\label{prop:explicit_transgression}
  We have $\tilde T(x)=T(x)$.
\end{proposition}
\begin{proof}
  The recipe described in the proposition defines a transformation
  $\tilde T$ which is again natural with respect to the map
  $f:E\rightarrow B$. As explained in the Remark
  \ref{ref:suspension_determines_transgression}
it must coincide with the transgression $T$ if it does so
in the special case
   of the cofibration $B\to
  CB\to\Sigma B$.  But in this case the
  above description produces
  exactly the suspension isomorphism which is by definition the
  transgression map $T$.
\end{proof}  

\subsubsection{}
If $f:E\rightarrow B$ is a map of smooth manifolds and $R=\R$,
then we could replace the singular cochains by differential forms
and construct $\tilde T$ on the level of forms. Again we get
$\tilde T=T$.

\subsubsection{} 
\label{sec:transgr-bockst}
For a cohomology class $x\in H^*(X,\integers)$ let $x^\rationals\in H^*(X,\rationals)$
  denote the image of $x$  under the canonical
  coefficient homomorphism $\integers\to\rationals$.
Let $n\in\nat $ and $\beta$ be the Bockstein transformation associated
to the sequence $$0\to \integers\xrightarrow{\cdot
    n}\integers\to \integers/n\integers\to 0\ .$$
 
\subsubsection{}
Let $x\in H^k(X,\integers)$ be such that
$nx$ is transgressive. Note that then $q^\Q$ is transgressive, too.
Since $nf^*x=0$ by the Bockstein exact sequence there exist 
    $u\in H^{k-1}(E,\integers/n\Z)$ with $\beta(u)=f^*x$. Recall
that $T(nx)$ is an equivalence class of cohomology classes. 
\begin{proposition}\label{prop:transgression and Bockstein} 
\begin{enumerate}
\item
We have
  \begin{equation}\label{eq:beta and transgression}
     T(nx) \ni i^* u ,
  \end{equation}
   
\item
     \begin{equation}\label{eq:rational beta and transgression}
    n T(x^\rationals) \ni  i^* u^\rationals.
  \end{equation}
\end{enumerate}
\end{proposition}
\begin{proof}
We use the description of the transgression  on the
singular cochain level given in Proposition \ref{prop:explicit_transgression}. 
Let $c$ be an
  integral cocycle representing $x$. Let $c_0$ be
  an integral cochain of $E$ with $dc_0 = nf^*c$. Then 
  $i^*c_0$ represents $T(nx)$. 

  The reduction of $c_0$ modulo $n\integers$ becomes closed and therefore
  represents a cohomology class $u\in
  H^{k-1}(E,\integers/n\integers)$. By the explicit construction of
  the Bockstein homomorphism, $\beta(u)=f^*x$. Equation \eqref{eq:beta and
    transgression} follows.

Since transgression commutes with the
  passage to rational coefficients by Lemma \ref{lem:transgression is
    compatible with natural transformations}, Equation
  \eqref{eq:rational beta and transgression} follows from Equation
  \eqref{eq:beta and transgression}.
\end{proof}

\section{Cohomology of $BO$, $BU$ and their loop spaces}
\label{sec:cohomology-bo-bu}
\subsection{The cohomology}
\subsubsection{}
In this appendix, we summarize the main results about the cohomology
of $BO$, $BU$ and their loop spaces, and the relations between them,
including the determination of the transgression homomorphisms. These
results are all classical, and almost all of them can be found in Cartan's
\cite{cartan}, where these calculations were
essential in his cohomological proof of Bott periodicity. Since they
are scattered over these papers, we collect them here in more
convenient form. All results without a proof or a different reference
can be found in \cite{cartan}.
\subsubsection{}
Bott periodicity gives canonical (up to homotopy) homotopy
  equivalences between $\Omega^n BO$ and other classical spaces
  summarized in the following list.
\begin{theorem}\label{theo:Bott periodicity list}
$$\begin{array}[l]{l|ccccccccccccc}
n-1 & -1 & 0 & 1 & 2 & 3 & 4 & 5 & 6 & 7  \\ \hline
  \Omega^nU/O & U/O & BO\times\integers & O & O/U & U/Sp &
  BSp\times\integers & Sp & Sp/U & U/O\\
  or &  & & SO\times\integers/2& SO/U\times\integers/2 
\end{array}$$
This extends $8$-periodically. 
\end{theorem}

  In the following, we will frequently identify the (cohomology of)
  different loop spaces of spaces in this table using the
  corresponding homotopy equivalence without further mentioning
  it. Note that we have done so already throughout the body of the paper.
 
\subsubsection{}
In the following, $L(x_{i_1},x_{i_2},\dots)$ denotes a polynomial algebra in
  the generators $x_i$, where by convention $x_i$ has (cohomological)
  degree $i$, and $E(y_{i_1},y_{i_2},\dots)$ denotes an exterior algebra, with similar
  degree conventions for the generators.
 
\subsubsection{}
 In the following list, we describe the cohomology of the connected
  component of the base point in $\Omega^kBO$. Note that we ``rename''
  some of the usual characteristic classes like the Pontryagin
  classes: $p_4$ is a cohomology class in $H^4$ etc.
\begin{theorem}\label{theo:cohom of BO etc}
$$\begin{array}{lcccl}
k & \Omega^kBO & H^*(\Omega^kBO_0,\integers) & H^*(\Omega^k BO_0,
\integers/2\Z) & H^*(\Omega^k BO_0,\integers[\frac{1}{2}]) \\ \hline
0 & BO &  L(p_{4},p_8,\dots)\oplus\text{$2$-Tors} & L(w_1,w_2,\dots) &
L(p_4,p_8,\dots) \\ 
1& O & * & L(d_1,d_3,\dots) & E(v_3,v_7,\dots) \\
2& O/U & L(u_2,u_6,\dots) &\\
3& U/Sp & E(a_1,a_5,\dots)\\
4& BSp\times\integers & L(y_4,y_8,\dots)\\
5& Sp & E(y_3,y_7,\dots)\\
6 & Sp/U & L(u_2,u_6,\dots)^* & E(c_2,c_4,\dots) & L(c_2,c_6,\dots)\\
6 &Sp/U & \displaystyle{\frac{L[c_2,c_4,\dots]}{\sum_{i+j=2k} (-1)^i c_{2i}c_{2j}}}\\
7 &U/O& E(a_{4k+1})\oplus\text{$2$-Tors} & E(w_1,w_2,\dots) & E(a_1,a_5,\dots)\\
\end{array}
$$
We add the following detailed explanations, using the description of
the loop spaces as in Theorem \ref{theo:Bott periodicity list}.
\begin{enumerate}
\item
  \begin{enumerate}
  \item 
    $H^*(BO;\integers)$ contains a subalgebra isomorphic to the
    quotient by its torsion.
  \item This is a polynomial algebra
    $L(p_4,p_8,\dots)$. 
  \item The torsion is annihilated by $2$, it is the
    image of Bockstein. 
  \item Reduction mod $2$ maps $p_{4k}$ to
    $(w_{2k})^2$. 
  \item The classes $w_{2k+1}\in H^{2k+1}(BO,\integers/2\Z)$
    have unique lifts to $H^{2k+1}(BO,\integers)$ which we also denote
    by $w_{2k+1}$. 
  \item The same is true for every class in degree $k$ for $k$
    not divisible by $4$, since in these degrees $H^k(BO,\integers)$ is
    annihilated by $2$.
\end{enumerate}
\item Most complicated is the cohomology of $SO$ with
  $\integers$-coefficients, for reasons of space simply denoted $*$ in
  the
  list (case $k=1$). We can say the following about it.
  \begin{enumerate}
  \item 
    The torsion in $H^*(SO,\integers)$ is annihilated by $2$, it is the
    image of Bockstein.
  \item The quotient of $H^*(SO,\integers)$ by its torsion is an
    exterior algebra $E(\overline{v_3},\overline{v_7},\dots)$. It does
    not split back to $H^*(SO,\integers)$ because of the product
    structure (compare with $H^*(SO,\integers/2\Z)$).
  \item  But of course, each
    monomial $\overline{v_{i_1}}\dots \overline{v_{i_s}}$ has an
    inverse image $v_{i_1}\dots v_{i_s}\in H^*(SO,\integers)$ (only
    additive! no multiplicative structure) which is well defined up to
    torsion, and the products are correct up to torsion.
\end{enumerate}
\item The integral cohomology of $Sp/U$ (case $k=6$) is the dual of
  $L(u_2,u_6,\dots)$. This shows in particular, that it is
  torsion-free. As a ring, it is the quotient of $L(c_2,c_4,\dots)$ by
  the ideal generated by the elements $\sum_{i+j=2k} (-1)^i c_{2i}c_{2j}$.
\item 
  \begin{enumerate}
  \item $H^*(U/O,\integers)$ contains a subalgebra isomorphic to the
    quotient by its torsion, this is an exterior algebra
    $E(a_1,a_5,\dots)$.
  \item The torsion is annihilated by $2$, it is the image of
    Bockstein.
  \item Reduction mod $2$ maps $a_{4k+1}$ to $w_{2k} w_{2k+1} +
    \beta(w_{4k} + w_2 w_{4k-2} +\dots+ w_{2k-2}w_{2k+2})$.
\end{enumerate}
\end{enumerate}
\end{theorem}

\subsubsection{}
 We also need the complex case, i.e.~$BU$ (and will later
  relate $BO$ to $BU$). The case of $BU$ is
of course much easier because of $2$-periodicity, and since the
cohomology does not contain torsion.
\begin{theorem}\label{theo:complex Bott}
$$\begin{array}{lccll}
  k & \mbox{natural homotopy equivalence of $\Omega^k BU$ to }& H^*(\Omega^k
  BU_0, \integers),\\ \hline 
  0 & BU & L(c_2,c_4,\dots)\\
  1 & U  & E(c_1,c_3,\dots)\\
  2 & BU\times\integers.
\end{array}$$
\end{theorem}

\subsubsection{}

  We now describe the effect of the loop map $\Omega\colon H^*(X)\to
  H^{*-1}(\Omega X)$ for integral cohomology and some of the spaces in \ref{theo:Bott
    periodicity list}. 

For the following table, recall that the universal Chern character is a
certain rational polynomial in the universal Chern classes, and we have a
unique integral lift $k! \ch_{2k}\in H^{2k}(BU\times\integers;\integers)$. 

\begin{theorem}\label{theo:loop map for BO etc}
$$\begin{array}{ccccl}
  \mbox{Space $X$} & \Omega X & x\in H^*(X,\Z) & \Omega(x)\in H^{*-1}(\Omega
  X,\Z)\\ \hline
  BU & U & c_{2k} & c_{2k-1}\\
  U & BU\times \Z & c_{2k-1} & (k-1)!\ch_{2k-2} \\
  BSp & Sp & y_{4k} & y_{4k-1}\\
  BO & O & p_{4k}& 2 v_{4k-1}+ \text{Tors}
\end{array}$$
\end{theorem}
\begin{proof}
  We only have to prove that $\Omega (c_{2k-1}) = (k-1)!\ch_{k-2}$. For this, observe
  that 
  \begin{equation*}
c_{2k-1}=\Omega( c_{2k})= \Omega ((k-1)!\ch_{2k}) =
  (k-1)!\Omega(\ch_{2k}),
\end{equation*}
since the other summands in $\ch_{2k}$ are decomposable
  and because the loop map applied to a decomposable class is zero by
  Proposition \ref{prop:transgression_of_products}. Now the Chern
  character is compatible with Bott periodicity, and therefore
  $\Omega^2(\ch_{2k})=\ch_{2k-2}$. Consequently
  \begin{equation*}
    \Omega c_{2k-1} = (k-1)! \Omega^2(\ch_{2k}) = (k-1)!\ch_{2k-2}.
  \end{equation*}
\end{proof}

\subsubsection{} 
\begin{lem}\label{lem:l factors through Bockstein}
The natural map 
\begin{equation*}
l_*:H^{4k+1}(BO,\Z/2\Z)\rightarrow
H^{4k+1}(BO,\Q/\Z)
\end{equation*}
of Equation \eqref{eq:def_of l} factors through the image
of the Bockstein homomorphism
$$\beta:H^{4k+1}(BO,\Z/2\Z)\rightarrow H^{4k+2}(BO,\Z).$$
\end{lem}
\begin{proof}
  We have the following map of long exact sequences
  $$\begin{array}{ccccccccc}
    \stackrel{2\times_*}{\rightarrow}&H^{4k+1}(BO,\Z)&\stackrel{q_*}{\rightarrow}& H^{4k+1}(BO,\Z/2\Z)&\stackrel{\beta}{\rightarrow}&H^{4k+2}(BO,\Z)&\rightarrow\\
    &\frac{1}{2}\times_*  \downarrow&&l_*\downarrow&&\|\\
    \rightarrow & H^{4k+1}(BO,\Q)&\rightarrow &H^{4k+1}(BO,\Q/\Z)&
    \stackrel{\beta^\prime}{\rightarrow} &H^{4k+2}(BO,\Z)&\rightarrow
\end{array}
\ .$$
The assertion now follows from the fact that
$H^{4k+1}(BO,\Q)=0$.
\end{proof}  

\subsection{Maps between loop spaces of $BO$}
\label{sec:maps-between-our}
\subsubsection{}
There is a large number of canonical maps between the different spaces
in Theorem \ref{theo:Bott periodicity list} and in Theorem
\ref{theo:complex Bott} which are important for us and which are
described in the following list.
\begin{enumerate}
\item\label{item:O_to_U} The inclusion $c\colon O\to U$ (given by complexification).
\item\label{item:BO_toBU} The induced map $Bc\colon BO\to BU$, which gives rise
  to $Bc\times \id_\integers\colon BO\times\integers\to BU\times \integers$
\item The inclusion $q\colon U\to Sp$ (given by tensoring with the
  quaternions)
\item the induced map $Bq\colon BU\to BSp$
\item the inclusion $f\colon U\to O$ (given by forgetting the complex
  structure)
\item the induced map $Bf\colon BU\to BO$
\item\label{item:Sp_to_U} the inclusion $j\colon Sp\to U$ (given by
  forgetting the quaternionic   structure)
\item\label{item:Bsp_toBU} the induced map $Bj\colon BSp \to BU$, which gives rise to
  $Bj\times\id_{\integers}\colon BSp\times\integers\to
  BU\times\integers$.
\item the projection $p\colon U\to U/O$
\item the projection $P\colon U\to U/Sp$
\item\label{item:U_modO_toBU} the inclusion of the fiber $i\colon
  U/O\to BO$ obtained by dividing the
  total space of the universal principal $U$-fibration $U\to EU\to BU$
  by $O$ (here we use the fact that $EU/O$ is a model for $BO$)
\item \label{item:projection BO to BU} The fibration $BO\to BU$
  constructed in \ref{item:U_modO_toBU}
\item the map $Omega^7_i\colon Sp/U\to U/O$ obtained by looping this fibration
  seven times, it 
  is the fiber of the map $U/O\to BO$ obtained by looping the fibration $BO\to
  BU$ seven times. 
\item\label{item:Sp_modU_to_BU} the (similar) inclusion of the fiber
  $I\colon Sp/U \to BU$.
\item\label{item:O_mod_U_to_BU} the (similar) inclusion of the fiber
  $\phi\colon O/U\to BU$.
\item\label{item:U_modO_to_U} a map $\alpha\colon U/O\to U$ given as composition 
  \begin{equation*}
U/O\xrightarrow{\sim}\Omega Sp/U\xrightarrow{\Omega I}
  \Omega BU \xrightarrow{\sim} U
\end{equation*}
where the first map is the Bott periodicity homotopy
  equivalence, and the third is the usual homotopy equivalence (which
  is also part of (complex) Bott periodicity).
\item\label{item:U_modSp_to_U} a similar map $\beta\colon U/Sp\to U$,
  given as composition
  \begin{equation*}
  U/Sp\xrightarrow{\sim} \Omega O/U\xrightarrow{\Omega \phi} \Omega BU \xrightarrow{\sim} U
\end{equation*}
\end{enumerate}
\subsubsection{}
The following relations hold between these maps. As usual, we will
freely use the Bott periodicity homotopy equivalences of Theorem
\ref{theo:Bott periodicity list} and Theorem \ref{theo:complex Bott}
to identify certain loop spaces with other spaces (therefore, strictly
speaking, the following assertions are true up to homotopy).
\begin{enumerate}
\item It is a general fact in the theory of classifying spaces that
  one way to construct $Bc$ in \ref{item:BO_toBU} is the fibration map of
  \ref{item:projection BO to BU}, which therefore can be identified
  with $Bc$. Reason: the identity map $EU\to EU$, where the domain is
  considered as contractible $O$-principle bundle and the target as
  contractible $U$-principle bundle intertwines, using the inclusion
  $c\colon O\to U$ the structures as principle bundles. Therefore the
  induced map on the quotients is the map $Bc$.
\item The map \ref{item:O_to_U} is obtained from \ref{item:BO_toBU} by
  applying the loop space transformation (and using the Bott
  periodicity identifications of $\Omega BO\xrightarrow{\sim} O$ and
  $\Omega BU\xrightarrow{\sim} U$).
\item  
 Similarly,
  \ref{item:Sp_to_U} is obtained by applying the loop space functor to
  \ref{item:Bsp_toBU}.
\item  By construction, looping
  \ref{item:Sp_modU_to_BU} gives \ref{item:U_modO_to_U}.
\item By construction, looping
  \ref{item:O_mod_U_to_BU} gives \ref{item:U_modSp_to_U}.
\item  
  Cartan \cite{cartan} proves that $O/U\to BU$, i.e.~~\ref{item:O_mod_U_to_BU} is
  obtained by applying the loop space functor to the inclusion \ref{item:O_to_U},
  $O\to U$.
\item Cartan \cite{cartan} proves that looping \ref{item:Sp_to_U} gives
  \ref{item:Sp_modU_to_BU}. This requires to check that his explicitly
  given maps $Sp/U\to BU$ and $O/U\to BU$ are the fiber inclusions we
  claim they are. 
\item Cartan \cite{cartan} also checks that looping \ref{item:U_modSp_to_U}
  gives \ref{item:Bsp_toBU}.
\item Cartan \cite{cartan} proves that looping \ref{item:U_modO_to_U}
  gives \ref{item:BO_toBU}. Strictly speaking, in this and the
  previous case he considers the
  corresponding maps of universal coverings $SU/Sp\to SU$ which loops
  to $BSp\to BU$, and $SO/SU\to SO$ which loops to $BO\to BU$. Since
  we know that $U/Sp\to U$ and $O/U\to O$ induce isomorphisms on
  $\pi_1$-level (all isomorphic to $\integers$), the claim follows.
\end{enumerate}

\subsubsection{}
To conclude, we have shown that in the sequence
\begin{align*}
Bc\times \id &\colon BO\times\integers \to BU\times\integers\\
c&\colon O\to U\\
\phi& \colon O/U\to BU\\
\beta &\colon U/Sp\to U\\
Bj\times \id &\colon BSp\times\integers\to BU\times\integers\\
j &\colon Sp\to U\\
I&\colon Sp/U\to BU\\
\alpha&\colon U/O\to U\\
Bc\times \id&\colon BO\times\integers\to BU\times\integers\\
\end{align*}
each map is obtained by looping the previous one (and applying Bott
periodicity to identify the loop spaces with the next spaces in the list).

\subsubsection{}
In the following table, we list the effect of the maps in
cohomology. Again, this is due to Cartan \cite{cartan}, with a few exceptions easily
obtained from his work. In these cases, the reason is indicated in the
last column of the following table. Recall that we always only
consider the cohomology of the connected component of the base
point. ``By looping'' means that we know that certain maps are obtained
from each other by applying the loop space functor (and some canonical
homotopy equivalences), and that we know the effect of the natural
loop map functor $\Omega\colon H^*(X) \to H^{*-1}(\Omega X)$ by
Theorem \ref{theo:loop map for BO etc}. 
 
\subsubsection{}
\mbox{}\newline
\noindent
 \begin{tabular}{lllp{7cm}}\label{table:maps-between-loop their effect in homology}
  $f\colon X\to Y$ & $x\in H^*(Y,\Z)$ & $f^*(x_0) \in H^*(X_0,\Z)$ & reason\\ \hline
  $BSp\to BU$ & $c_{4k} $ & $y_{4k}$\\
  & $c_{4k+2}$ & $0$\\
  $O/U \to BU$ & $c_{4k+2}$ & $2u_{4k+2}$\\
               & $c_{4k}$ & \multicolumn{2}{l}{$-2\cdot \sum_{0<i<2k}
                 (-1)^i f^*(c_{2i})/2\cdot f^*(c_{4k-2i})/2 $} \\
               $Sp/U \to BU$ & $c_{2k}$ & $\overline{c_{2k}}$\\
               &      $c_{2k} \mod 2$ & $c_{2k}$\\
  $Sp \to U$ & $c_{4k+1}$ & $0$ & by looping\\
    & $c_{4k+3}$ & $ y_{4k+3}$ & by looping\\
  $U/Sp \to U$ & $c_{4k+1}$ & $2a_{4k+1}$&  (since true dually in
  homology)\\
               & $c_{4k+3}$ & $0$  & (since true dually in homology)\\
  $U\to Sp$ & $y_{4k-1}$ & $2c_{4k-1}$ & by looping (since products
  suspend to zero)\\
  $BU\to BSp$ & $y_{4k}$ & $\sum_{i+j =2k} (-1)^i c_{2i} c_{2j}$\\
  $SO\to U$  & $c_{4k+3}$ & $2{v_{4k+3}}+\text{$2$-Tor}$\\
             & $c_{4k+1}$ & $\text{$2$-Tor}$\\
  $U\to U/O$ & $a_{4k+1}$ & $2c_{4k+1}$\\
  $U/O\to U$ & $c_{4k+1}$ & $a_{4k+1}+\text{Tors}$\\
             & $c_{4k+3}$ & Tors\\
  $U\to O$ & ${v_{4k+3}} $ & $c_{4k+3}$\\
  $BO\to BU$ & $c_{4k}$ & $p_{4k}$\\
  & $c_{4k+2}$ & $w_{2k+1}^2$\\
  & $c_{4k}$ mod $2$ & $w_{2k}^2$\\
  $BU\to BO$ & $p_{4k}$ & $\sum_{i+j=k} (-1)^i c_{2i}c_{2j}$\\
  &  $w_{2k+1}$ & $0$\\
  &  $w_{2k}$ & $c_{2k}$ mod $2$\\
  $U/O\to BO$ & $w_k$ mod $2$ & $w_k$ mod $2$ \\
              & $p_{4k}$ & $0$ &  mod $2$ it maps to $0$, and pulled
              back further to $U$ it is also $0$, i.e.~no $2$-torsion
              and no free part\\
\end{tabular}
 
\subsubsection{}
In two cases, we
 have to take the different components into account: note that $\ch_0\in
 H^0(BU\times \integers;\integers)$ has the value $d$ (times the canonical
 generator) on the component of $BU\times\integers$ labeled by
 $d\in\integers$. Correspondingly, we have a real version $\ch_0^{\reals}\in
 H^0(BO\times\integers;\integers)$ and a quaternionic version
 $\ch_0^{\quaternions} \in H^0(Bsp\times\integers;\integers)$, describing the
 dimension of the virtual universal real or quaternionic bundle,
 respectively. For these classes we get the
 (obvious) following relations under the maps induced by ``complexification''
 or ``forgetting the  quaternionic structure'', respectively.

\begin{tabular}{lllllll}
  $f\colon X\to Y$ & $x\in H^0(Y,\Q)$ & $f^*(x) \in H^0(X,\Q)$\\ \hline
  $BSp\times\integers\to BU\times\integers$ & $\ch_{0}$ & $2\ch^{\quaternions}_{[0]}$\\
  $BO\times\integers\to BU\times\integers$ & $\ch_{0}$ & $\ch^{\reals}_{[0]}$\\
\end{tabular}
\subsubsection{}
\begin{lemma}\label{lem:transgressive classes}
  In the fibration $U/O\xrightarrow{i} BO \xrightarrow{Bc} BU$, the
  classes $c_{4k+2}^{\rationals}\in H^{4k+2}(BU;\rationals)$ are transgressive.
\end{lemma}
\begin{proof}
  The pull back of $c_{4k+2}$ to $BO$ is $2$-torsion in integral
  cohomology, therefore vanishes in rational cohomology.
  \end{proof}  

\begin{lemma}\label{lem:surjecti}
   In the fibration $Sp/U\xrightarrow{\Omega^7i} U/O \xrightarrow{\alpha} U$,
   the Leray-Serre spectral sequence for $H^*(\cdot,\integers/2\integers)$
   collapses at the $E_2$-term. In particular, the edge homomorphism
   \begin{equation*}
   (\Omega^7 i)^* \colon H^*(U/O;\integers/2\integers)\to
   H^*(Sp/U;\integers/2\integers)
 \end{equation*}
is surjective, whereas the edge
   homomorphism
   \begin{equation*}
\alpha^*\colon H^*(U;\integers/2\integers)\into
   H^*(U/O;\integers/2\integers)
 \end{equation*}
is injective.
\end{lemma}
\begin{proof}
  $E_2 = H^*(U/O;\integers/2\integers)\tensor H^*(Sp/U;\integers/2\integers)$
  is the tensor product of an exterior algebra over $\integers/2\integers$
  with exactly one generator in
  each positive even degree with an exterior algebra over $\integers/2$ with
  one generator in each positive odd 
  degree, i.e.~an exterior algebra over $\integers/2\integers$ with one
  generator in each positive degree. It converges to an exterior algebra over
  $\integers/2\integers$ with one generator in each positive degree. Any non-zero differential would
  result in an $E_\infty$-term which has is too low dimensional, therefore the
  spectral sequence necessarily collapses. The statement about the edge
  homomorphisms is an immediate consequence.
\end{proof}

\end{appendix}

\bibliographystyle{plain}
%\bibliography{/home/uli/working/libank/literatu,/home/uli/working/libank/lit1} 

\begin{thebibliography}{10}

%\bibitem{atiyah???}
%M.~F. Atiyah.
%\newblock {\em $K$-theory} 
%\newblock W.A. Benjamin, Inc., New York-Amsterdam, 1967.

\bibitem{atiyahhirzebruch}
M.~F. Atiyah and F. Hirzebruch.
\newblock Vector bundles and homogeneous spaces.
\newblock {\em Proceedings of Symp. in pure Math}, Vol.3, 1961.

\bibitem{atiyahsinger69}
M.~F. Atiyah, V.~K. Patodi, and I.~M. Singer.
\newblock Index theory for skew-adjoint Fredholm operators.
\newblock {\em Publ. I.H.E.S}, Vol. 37 (1969), 5--26.

\bibitem{atiyahpatodisinger75}
M.~F. Atiyah, V.~K. Patodi, and I.~M. Singer.
\newblock Spectral asymmetry and Riemannian geometry I.
\newblock {\em Math.Proc.Camb.Phil.Soc.}, 77(1975), 43--69, 1975.

\bibitem{berlinegetzlervergne92}
N.~Berline, E.~Getzler, and M.~Vergne.
\newblock {\em Heat Kernels and Dirac Operators}.
\newblock Springer-Verlag Berlin Heidelberg New York, 1992.

%\bibitem{bismut86}
%J.~M. Bismut.
%\newblock The index theorem for families of Dirac operators: Two heat equation
%  proofs.
%\newblock {\em Invent. Math}, 83(1986), 91--151.

\bibitem{bismutcheeger89}
J.~M. Bismut and J.~Cheeger.
\newblock $\eta$-invariants and their adiabatic limits.
\newblock {\em J. AMS}, 2(1989), 33--70.

%\bibitem{bismutcheeger90}
%J.~M. Bismut and J.~Cheeger.
%\newblock Families index for manifolds with boundary, superconnections, and
%  cones. I. Families of manifolds with boundary and Dirac operators.
%\newblock {\em J. Funct. Anal.}, 89(1990), 313--363.

%\bibitem{bismutcheeger901}
%J.~M. Bismut and J.~Cheeger.
%\newblock Families index for manifolds with boundary, superconnections, and
%  cones. II. The Chern character.
%\newblock {\em J. Funct. Anal.}, 90(1990), 306--354.

%\bibitem{bismutfreed861}
%J.~M. Bismut and D.~Freed.
%\newblock The analysis of elliptic families. I. Metrics and connections on
%  determinant bundles.
%\newblock {\em Comm. Math. Phys.}, 106(1986), 159--176.

%\bibitem{bismutfreed862}
%J.~M. Bismut and D.~Freed.
%\newblock The analysis of elliptic families. II. Dirac operators, eta
%  invariants, and the holonomy theorem.
%\newblock {\em Comm. Math. Phys.}, 107(1986), 103--163.

%\bibitem{borelhirzebruch59}
%A. Borel and F. Hirzebruch.
%\newblock Characteristic classes and homogeneous spaces II.
%\newblock {\em Amer. J. of Math.}, 81 (1959), 315--381.

\bibitem{brylinski93}
J.~L. Brylinski.
\newblock {\em Loop spaces, characteristic classes, and geometric
  quantization}.
\newblock Birkh{\"a}user, Progress in Math. 107, 1993.

\bibitem{bunke020}
U.~Bunke.
\newblock $\eta$-form, index theory and Deligne cohomology.
\newblock Preprint 2002. arXiv:math.DG/0201112

%\bibitem{bunke011}
%U.~Bunke.
%\newblock Transgression of the index gerbe.
%\newblock {\em manuscripta math.}, 109 (2002), 263-287. arXiv:math.DG/0109052

%\bibitem{bunkekoch98}
%U.~Bunke and H.~Koch.
%\newblock The $\eta$-form and a generalized Maslov index.
%\newblock {\em manuscripta math}, 95(1998), 189--212. arXiv:dg-ga/9701004

\bibitem{bunkema02}
U.~Bunke and X.~Ma.
\newblock Index and secondary index theory for flat bundles with duality.
\newblock Preprint 2002 : http://www.uni-math.gwdg.de/bunke/sample3.pdf

\bibitem{cartan}
H.~Cartan.
\newblock{D{\'e}monstration homologique des th{\'e}or{\`e}mes de p{\'e}riodicit{\'e} de Bott. I, II, III}
\newblock{P{\'e}riodicit{\'e} des Groupes d'Homotopie stables des Groupes
  classiques, d'apres Bott, {\em S{\'e}m. H. Cartan}  12 (1959/60), No.16, 16 p., No.17, 32 p., No.18, 9 p. }
({1961})

\bibitem{cheegersimons83}
J.~Cheeger and J.~Simons.
\newblock Differential characters and geometric invariants.
\newblock In {\em LNM1167}, pages 50--80. Springer Verlag, 1985.

%\bibitem{daizhang???}
%X.~Dai and W.~Zhang.
%\newblock Higher spectral flow.
%\newblock {\em J. Funct. Anal.}, 157 (1998), 432-469. arXiv:dg-ga/9608002 

\bibitem{freed99}
D.~Freed.
\newblock Dirac charge quantization and generalized differential cohomology.
\newblock {\em In: Surveys in differential geometry}, International
Press, Somerville, 2000, 129-194. arXiv:hep-th/0011220

%\bibitem{gajer97}
%P.~Gajer.
%\newblock Geometry of Deligne cohomology.
%\newblock {\em Invent. Math.}, 127(1997), 155--207. arXiv:alg-geom/960125

%\bibitem{gajer99}
%P.~Gajer.
%\newblock Higher holonomies, geometric loop groups and smooth Deligne
%  cohomology.
%\newblock In {\em Advances in geometry. Progr. in Math. 172}, pages 195--235.
%  Birkh{\"a}user, Boston, 1999.

%\bibitem{goette00}
%\newblock Equivariant $\eta$-invariants and $\eta$-forms.
%\newblock {\em J. reine angew. Math.} 526(2000), 181--236. arXiv:math.DG/0203269

%\bibitem{hiltonwylie60}
%P.~J. Hilton and S. Wylie
%\newblock{\em Homology theory.}
%\newblock Cambridge Univ. Press, 1960.

\bibitem{hitchin74}
N.~Hitchin.
\newblock Harmonic spinors.
\newblock {\em Adv. in Math.} 14 (1974), 1-55.


%\bibitem{hitchin99}
%N.~Hitchin.
%\newblock Lectures on special Lagrangian submanifolds.
%\newblock arXiv:math.DG/9907034

%\bibitem{lott941}
%J.~Lott.
%\newblock {$\R/\Z$}-index theory.
%\newblock {\em Comm. Anal. Geom.}, 2(1994), 279--311.

%\bibitem{lott01}
%J.~Lott.
%\newblock Higher degree analogs of the determinant line bundle.
%\newblock {\em Comm. Math. Phys.}, 230(2002), 41-69. arXiv:math.DG/0106177

\bibitem{karoubi78}
M.~Karoubi.
\newblock $K$-theory.
\newblock Springer-Verlag, 1978.


%\bibitem{kervaire59}
%M.~A. Kervaire.
%\newblock A note on obstructions and characteristic classes.
%\newblock {\em Amer. J. of Math.}, 81(1959), 773--784.

%\bibitem{melrose93}
%R.~B. Melrose.
%\newblock {\em The Atiyah-Patodi-Singer Index Theorem}.
%\newblock A.K.Peters, Wellesley, 1993.

%\bibitem{melrosepiazza97}
%R.~B. Melrose and P.~Piazza.
%\newblock Families of Dirac operators, boundaries, and the $b$-calculus.
%\newblock {\em J. Differential. Geom.}, 46(1997), 99--180.

\bibitem{milnorstasheff68}
J.W.~Milnor and J.D.~Stasheff.
\newblock {\em Characteristic classes.}
\newblock Annals of Mathematics Studies, Vol 76, Princeton University
Press, Princeton, 1974.

\end{thebibliography}
%bibliography{/user/bunke/working/libank/literatu,/user/bunke/working/libank/lit1}

\end{document}